%% file: Balabanov_Beaupere_Grigori_Lederer_blockSRHT.tex
\documentclass[10pt]{article}
\usepackage[includeheadfoot,top=2.97 cm, bottom=2.97 cm, left=2.1 cm, right=2.1 cm]{geometry}
\usepackage{appendix}

\usepackage[utf8]{inputenc} 
\usepackage[T1]{fontenc}    
\usepackage{hyperref}      
\usepackage{url}           
\usepackage{booktabs}       
\usepackage{amsfonts}       
\usepackage{nicefrac}       
\usepackage{microtype}      
\usepackage{xcolor}        

\usepackage{tikz}
\usetikzlibrary{datavisualization}
\usetikzlibrary{plotmarks}
\pgfdvdeclarestylesheet{cust}
{
	default style/.style={visualizer color=black},
	1/.style={visualizer color=red},
	2/.style={visualizer color=black},
	3/.style={visualizer color=red},
	4/.style={visualizer color=black}
}
\usepackage{caption}
\usepackage{subcaption}

\usepackage{graphicx}
\usepackage{lipsum}
\usepackage[utf8]{inputenc}

\usepackage{algorithm}
\usepackage{algorithmicx}
\usepackage{algpseudocode}

\usepackage{amsthm}
\usepackage{upgreek}
\usepackage{amssymb}
\usepackage{hyperref}
\usepackage{multirow}
\usepackage{mathtools}
\usepackage{color}
\usepackage{xcolor}
\usepackage{float}
\usepackage{subcaption}
\usepackage[noabbrev,capitalise]{cleveref}

\Crefname{proposition}{Proposition}{Propositions}

\newtheorem{lemma}{Lemma}[section]
\newtheorem{proposition}[lemma]{Proposition}

\newtheorem{theorem}[lemma]{Theorem}
\newtheorem{remark}[lemma]{Remark}
\newtheorem{definition}[lemma]{Definition}
\crefname{equation}{}{}
\Crefname{equation}{}{}
\numberwithin{equation}{section}

\newcommand{\bphi}{\mathbf{\upvarphi}}

\newcommand{{{\bxi}}}{\mathbf{\upxi}}
\newcommand{\Frob}{\mathrm{F}}

\newcommand{\bb}{\mathbf{b}}

\newcommand{\by}{\mathbf{y}}

\newcommand{\bm}{\mathbf{m}}

\newcommand{\bz}{\mathbf{z}}
\newcommand{\bx}{\mathbf{x}}

\newcommand{\bTheta}{\mathbf{\Theta}}
\newcommand{\bSigma}{\mathbf{\Sigma}}

\newcommand{\bPhi}{\mathbf{\Phi}}
\newcommand{\bR}{\mathbf{R}}
\newcommand{\bA}{\mathbf{A}}

\newcommand{\bB}{\mathbf{B}}

\newcommand{\bhR}{\widehat{\mathbf{R}}}
\newcommand{\bhW}{\widehat{\mathbf{W}}}

\newcommand{\bhD}{\widehat{\mathbf{D}}}

\newcommand{\bhH}{\widehat{\mathbf{H}}}

\newcommand{\bX}{\mathbf{X}}
\newcommand{\bZ}{\mathbf{Z}}

\newcommand{\bM}{\mathbf{M}}

\newcommand{\bH}{\mathbf{H}}
\newcommand{\bD}{\mathbf{D}}

\newcommand{\bW}{\mathbf{W}}
\newcommand{\bOmega}{\mathbf{\Omega}}
\newcommand{\bGamma}{\mathbf{\Gamma}}

\newcommand{\bY}{\mathbf{Y}}
\newcommand{\bQ}{\mathbf{Q}}
\newcommand{\bC}{\mathbf{C}}
\newcommand{\bP}{\mathbf{P}}

\newcommand{\bU}{\mathbf{U}}

\newcommand{\bV}{\mathbf{V}}

\newcommand{\bPi}{\mathbf{\Pi}}

\newcommand{\bd}{\mathbf{d}}
\newcommand{\bPsi}{\mathbf{\Psi}}

\hypersetup{
	colorlinks,
	linkcolor={blue!},
	citecolor={blue!},
	urlcolor={blue!}
}
\def\BibTeX{{\rm B\kern-.05em{\sc i\kern-.025em b}\kern-.08em
		T\kern-.1667em\lower.7ex\hbox{E}\kern-.125emX}}
\begin{document}

	\title{Block subsampled randomized Hadamard transform for low-rank approximation on distributed architectures\\
	}
	
	\author{Oleg Balabanov, Matthias Beaup\`ere,	Laura Grigori and Victor Lederer \thanks{Sorbonne Universit\'e, Inria, CNRS, Universit\'e de Paris, Laboratoire Jacques-Louis Lions, Paris, France.}}
	\date{}
	\maketitle
	
	\begin{abstract}
		This article introduces a novel structured random matrix composed blockwise from subsampled randomized Hadamard transforms (SRHTs). The block SRHT is expected to outperform well-known dimension reduction maps, including SRHT and Gaussian matrices, on distributed architectures with not too many cores compared to the dimension. We prove that a block SRHT with enough rows is an oblivious subspace embedding, i.e., an approximate isometry for an arbitrary low-dimensional subspace with high probability. Our estimate of the required number of rows is similar to that of the standard SRHT. This suggests that the two transforms should provide the same accuracy of approximation in the algorithms. The block SRHT can be readily incorporated into randomized methods, for instance to  compute a low-rank approximation of a large-scale matrix. For completeness, we revisit some common randomized approaches for this problem such as Randomized Singular Value Decomposition and Nystr\"{o}m approximation, with a discussion of their accuracy and implementation on distributed architectures. 
	\end{abstract}
	
	
	\section{Introduction} \label{Intro}
	
	Randomization has become a powerful tool for tackling massive problems in numerical algebra and data science~\cite{martinsson2020randomized,woodruff2014sketching,vershynin2018high,mahoney2011randomized}. Modern randomized methods can, in particular, provide solutions to problems of dimensions beyond the reach of deterministic methods, and allow effective use of computational resources. Recent significant development has made them very reliable, and not just used as a last resort, as it was not so long ago. Along with increased efficiency, they can now provide strong accuracy guarantees with a user-specified failure probability that can be chosen extremely low, say $10^{-10}$, without much impact on computational costs.
	
	This article is concerned with randomized methods that are based on a dimension reduction, called sketching~\cite{woodruff2014sketching}, with oblivious $\ell_2$-subspace embeddings (OSEs) defined below. 
	\begin{definition} \label{def:oblemb}
		Let $0 \leq \varepsilon<1$ and $0< \delta <1$. A random matrix $\bOmega \in \mathbb{R}^{l \times n}$ is said to be a $(\varepsilon,\delta,d)$ OSE, if for any fixed $d$-dimensional subspace $V \subseteq \mathbb{R}^n$,
		\begin{equation} \label{eq:isometry}
		\forall \bx \in V,~~ | \| \bx \|_2^2 - \|\bOmega \bx\|_2^2 | \leq \varepsilon \| \bx \|_2^2
		\end{equation}
		holds with probability at least $1-\delta$.
	\end{definition}
	It is a consequence of the Johnson-Lindenstrauss lemma~\cite{johnson1984extensions} that there exist $(\varepsilon, \delta, d)$ OSEs of sizes $l = \mathcal{O}(\varepsilon^{-2} (d +
	\log{\frac{1}{\delta}})$. The fact that $n$ does not appear in the right-hand-side and the logarithmic dependence on the probability of failure $\delta$ shows the potential of a dimension reduction with such embeddings. There are several distributions that are known to satisfy the OSE property with the optimal or close to optimal $l$. The Gaussian, Rademacher distributions, sub-sampled randomized Hadamard transform (SRHT),  sub-sampled randomized Fourier transform, and CountSketch matrix are ones of the most popular distributions. The random sketching matrix in the algorithm should be chosen depending on the computational architecture to yield the most benefit. For instance, the SRHT is a structured matrix that can be efficiently applied to a vector in a sequential computational environment, while the application of a Rademacher matrix is efficient in a highly parallel environment. In this paper we propose a novel OSE, called block SRHT, which should be superior to all currently existing ones on a distributed computational architecture, with not too many computational units. 
	
	The OSEs are used in a variety of randomized methods for machine learning, scientific computing, and signal processing. Perhaps one of the most representative applications is the linear regression problem. Suppose that we seek a vector $\bx \in  \mathbb{R}^d$ 
	that minimizes $\| \bA \bx - \bb \|_2,$
	where $\bA \in \mathbb{R}^{n\times d}$ is a large-scale dense matrix, $\bb \in \mathbb{R}^n$ is a large-scale vector, and $d \ll n$.
	It follows that the solution to this problem can be approximated by a minimizer of $\| \bOmega (\bA \bx - \bb) \|_2$  requiring considerably lower computational cost.  The accuracy of such an approximation is guaranteed, given that $\bOmega$ is $(\varepsilon,\delta,d+1)$ OSE. Besides the linear regression problem, the sketching technique with OSEs has been successfully applied to the nearest neighbors problem~\cite{ailon2006approximate}, approximation of products of matrices~\cite{sarlos2006improved}, computation of low-rank approximations of matrices~\cite{halko2011finding} and tensor decompositions~\cite{sun2020low}, dictionary learning~\cite{anaraki2013compressive}, solution of parametric equations~\cite{balabanov2019randomized}, and solution of linear systems and eigenvalue problems~\cite{balabanov2020randomized}.

	In this paper the potential of the block SRHT is realized on the low-rank approximation problem. Such problems are ubiquitous, for instance, in the principal component analysis of large data sets and kernel ridge regression. A randomized low-rank approximation for machine learning tasks was addressed in e.g.~\cite{alaoui2015fast,bach2013sharp,derezinski2020improved}.  In~\cite{meanti2020kernel,yin2021distributed,rudi2017falkon,zhang2013divide,calandriello2016analysis} a particular focus was given to make the methods suited to modern architectures. In details, given a large matrix $\bA \in \mathbb{R}^{m \times n}$ with rapidly decaying spectrum, we seek a matrix $\bB_k$ preferably in an SVD form, of rank at most $k \ll \min(m,n)$, that approximates well $\bA$. The matrix $\bB_k$ can be obtained by first restricting its range to a subspace that captures the most of the action of $\bA$, and then minimizing the chosen error measure, say the spectral error $\|\bA-\bB_k\|_2$. As shown in~\cite{halko2011finding}, the most of $\bA$'s action can be well captured by the range of $\bA \bOmega^\mathrm{T}$, which constitutes the core of state-of-the-art Randomized Singular Value Decomposition (RSVD) algorithm. 
    Over the past years more sophisticated randomized low-rank approximation methods have been developed such as the Nystr\"{o}m method for spd matrices~\cite{tropp2017fixed,gittens2013revisiting}, and the single-view approximations for general matrices~\cite{upadhyay2016fast,tropp2019streaming,kannan2014principal}. In this work, we outline some such methods and analyze them with a projection-based approach compatible with block SRHT.
	
	The paper is organized as follows. The rest of~\Cref{Intro} discusses contributions and relation to prior work.  \Cref{blockSRHT} introduces a block SRHT matrix and discusses its properties. In~\cref{randlowrank} we present randomized algorithms for computing a low-rank approximation based on oblivious embeddings. \Cref{experiments} contains some computational aspects and experimental results. The proof of the main theoretical result is given in~\cref{mainproof}.  \Cref{conclusion} concludes this work. 
	
	\color{black}
	\subsection*{Contributions} \label{contributions}
	The proposed block SRHT matrix has the potential to combine the benefits of structured oblivious embeddings, such as the SRHT, with the benefits of unstructured ones, such as Gaussian, from the complexity and performance standpoint. We are not aware of any other work in this direction. We prove that the block SRHT matrix of size $l = \mathcal{O}(\varepsilon^{-2}(d+\log{\frac{n}{\delta}}) \log{\frac{d}{\delta}})$ satisfies the $(\varepsilon, \delta, d)$ OSEs property. This result is similar to that for standard SRHT from the literature~\cite{tropp2011improved,boutsidis2013improved} and implies that the two matrices should yield approximations of similar quality. This result does not simply follow from the analysis in~\cite{tropp2011improved,boutsidis2013improved}, and particularly requires incorporation of a useful technical trick that, to the best of our knowledge, was not employed before in the randomized numerical linear algebra community.

	In addition, we present low-rank approximation methods such as RSVD~\cite{halko2011finding}, Nystr\"{o}m approximation, and the single-view approximation from~\cite{tropp2019streaming} in a new unified projection-based form that clearly shows their connection. 
	We provide a rigorous characterization of their accuracy based solely on the OSE property, and therefore compatible with any types of sketching matrices that satisfy this property, including the block SRHT. For RSVD and Nystr\"{o}m approximation, this characterization follows almost directly from standard results from the literature. On the other hand, for the single-view approximation, our results are new. In particular, the novelty lies in the use of the projection-based interpretation of the single-view approximation to show that it is almost as accurate as RSVD if the embeddings for the ``core sketch'' are OSEs of sufficiently large sizes. Important aspects of implementation in distributed architectures are also discussed using the suitability of block SRHT for these architectures.
	
	\color{black}

	\section{Block sub-sampled randomized Hadamard transform} \label{blockSRHT}
	For $n$ being a power of two, an SRHT matrix can be defined as follows:
	\begin{equation} \label{eq:SRHT}
	\bOmega = \sqrt{\textstyle \frac{n}{l}}  \bR \bH \bD, 
	\end{equation}
	where $\bR$ is a $l \times n$ uniform, with or without replacement, random sampling matrix, $\bH$ is a $n \times n$ Walsh-Hadamard matrix rescaled by $\frac{1}{\sqrt{n}}$, and $\bD$ is a diagonal matrix with i.i.d. Rademacher random variables $\pm 1$ on the diagonal. 
	The properties of SRHT  were thoroughly described in~\cite{tropp2011improved} with a follow up analysis in~\cite{boutsidis2013improved}. SRHT matrices are commonly used in randomized algorithms as they can be applied to vectors using only $n \log_2n$ flops, while general unstructured matrices require $2 n l$ flops. At the same time, they satisfy the $(\varepsilon,\delta,d)$ OSE property, if~\cite{balabanov2019randomized}
	\small 
	\begin{equation} \label{eq:SRHTbound}
	l \geq 3 \varepsilon^{-2}  (\sqrt{d}+ \sqrt{8 \log{\textstyle \frac{6n}{\delta}}}  )^2 {\textstyle \frac{3d}{\delta}}, 
	\end{equation}
	\normalsize
	which is only by a logarithmic factor in $\delta$ and $n$ larger than the optimal bound.  For a general $n$, a partial SRHT can be used that is defined as the first $n$ columns of an SRHT matrix. Using a partial SRHT is equivalent to padding the input data with zeros to make its dimension a power of two. 
	
	{Unfortunately, products with SRHT matrices are not well suited to distributed computing limiting the benefits of SRHT on modern architectures {(see e.g.~\cite{yang2015implementing}).} This happens majorly due to computing products with $\bH$ in tensor form \small $$ \bH = {\textstyle \frac{1}{\sqrt{n}}}\begin{bmatrix} 1 &1 \\ 1 &-1 \end{bmatrix} \otimes \begin{bmatrix} 1 &1 \\ 1 &-1 \end{bmatrix} \otimes \hdots \otimes \begin{bmatrix} 1 &1 \\ 1 &-1 \end{bmatrix},$$
		\normalsize requiring cumbersome reduction operator such as a sequence of arrays of butterflies, rather than a simple addition, which we have with Gaussian matrices. This article attempts to alleviate this problem by constructing $\bOmega$ block-wise as follows
		\begin{equation} \label{eq:blockSRHT}
		\bOmega =  [ \begin{matrix}  \bOmega^{(1)}  & \bOmega^{(2)} & \hdots &  \bOmega^{(p)} \end{matrix} ],
		\end{equation} 
		where $\bOmega^{(i)} =\sqrt{\frac{r}{l}} \widetilde{\bD}^{(i)}   \bR \bH \bD^{(i)} $ are $l\times r$ SRHT matrices related to a unique sampling matrix $\bR$ and different (independent from each other) diagonal matrices $\bD^{(i)}$ with i.i.d. Rademacher entries $\pm1$, multiplied from the left by another diagonal matrices $\widetilde{\bD}^{(i)}$ with Rademacher entries, $r = \frac{n}{p}$, $1 \leq i \leq p$. As in the standard SRHT, the condition that $r$ is a power of two can be achieved by zero padding of the input data. 
		The advantage of $\bOmega$ defined by~\cref{eq:blockSRHT} is that it can be multiplied by an $n \times d$ matrix $\bV$  distributed between $p$ processors with rowwise partitioning, as 
		\begin{equation} \label{eq:blockSRHTappl}
		\bOmega \bV =   \sum_{1 \leq i \leq p} \bOmega^{(i)} \bV^{(i)},
		\end{equation}
		where $\bV^{(i)}$ are the corresponding local blocks of rows of $\bV$. In this way, to obtain $\bOmega \bV$ one can compute the local contributions $\bOmega^{(i)} \bV^{(i)}$ on each processor and then sum-reduce them to the master processor. {This makes block SRHT matrices have the same application cost in terms of communication as Gaussian matrices. Thus, they should yield much better scalability of computations than standard SRHT~\cite{yang2015implementing}.  The sum-reduce operation requires exchanging $\mathcal{O}(\log p)$ messages and $\mathcal{O}(dl \log p)$ per-processor communication volume that can be by a factor $\mathcal{O}(\frac{r}{l})$ less than the volume of communication used by standard SRHT (if $l \leq r$).}  {At the same time, block SRHT require less flops per processor than Gaussian matrices. To be more specific, the application cost of block SRHT using~\cref{eq:blockSRHTappl} is only $\mathcal{O}(rd \log{r}+ dl \log{p})$ flops per processor, while Gaussian matrices require $\mathcal{O}(rdl+ dl \log{p})$ flops per processor.} {It is deduced that block SRHT matrices are both well-suited to distributed computing and efficient in terms of flops. They are  expected to outperform all existing oblivious embedding when the local dimension $r$ and the sampling dimension $l$ are large enough.}

		The procedure for application of the block SRHT can be easily extended to the case when $\bV$ is distributed with a 2D partitioning. Namely, to multiply $\bOmega$ by $n \times n$ matrix $\bV$ distributed over a grid of $p \times p$ processors, we first compute the local contributions $\bX^{(i,j)}=\bOmega^{(i)} \bV^{(i,j)}$ on each processor, and then sum-reduce the contributions from the $j$-th column of blocks to the processor $(1,j)$, $1 \leq j \leq p$. Note that in this case the resulting matrix $ \bY^\mathrm{T} = \bOmega \bA$ is distributed with rowwise partitioning over processors $(1,1), (1,2), \hdots, (1,p)$. This provides the ability to efficiently compute the sketch $\bOmega \bY$, or to orthogonalize $\bY$ with a routine suited for distributed computing, with no need to reorganize $\bY$. This can be particularly handy in the low-rank approximation algorithm from~\cref{randlowrank}.


		We will assume that $\bR$ in~\cref{eq:blockSRHT} samples rows uniformly at random and  \emph{with replacement}.  Interestingly, in this case the block SRHT can be viewed  as a generalization of the SRHT with replacement and the Rademacher embedding, as it reduces to these maps when $r=n$ and $r=1$, respectively. Sampling with replacement can be important, for instance, when $r$ is smaller than the dimension of the embedded subspace.
		
		\Cref{thm:main} is the main result of the article. It implies the compatibility of the block SRHT with all randomized methods that rely on OSEs, including the methods in~\cref{randlowrank}.   The estimate of $l$ in~\cref{thm:main} is similar to~\cref{eq:SRHTbound} for the standard SRHT matrix, and in particular depends only logarithmically on {$n$} and $\delta$.
		\begin{theorem}[Main Theorem]\label{thm:main}
			Let $0< \varepsilon <1$ and $0< \delta <1$. Let $\bOmega \in \mathbb{R}^{l \times n}$ be defined by~\cref{eq:blockSRHT}.  If, \small $$ {n \geq l} \geq 3.7 \varepsilon^{-2} (\sqrt{d}+4\sqrt{\log{\textstyle \frac{n}{\delta}}+6.3})^2 \log{\textstyle \frac{5d}{\delta}},$$ \normalsize then $\bOmega$ is an $(\varepsilon,\delta,d)$ OSE.	
		\end{theorem}
		For better presentation the proof of~\cref{thm:main} is deferred to the end of the article (see~\cref{mainproof}).

		\section{Randomized low-rank approximation} \label{randlowrank}

        This section addresses the problem of efficient  computation of a rank-$k$ approximation of a large matrix $\bA\in \mathbb{R}^{m \times n}$ with rapidly decaying spectrum.  We provide and analyze some common randomized algorithms for this task. They are presented with a {$\emph{projection-based}$} approach relying solely on the OSE property of the sketching matrix, and are compatible with the block SRHT thanks to~\cref{thm:main}.
        {A particular focus is given to the scenario where $\bA$ is uniformly distributed over a 2D grid of processors. For simplicity assume that $m \leq n$.}

It is a well-known fact that the best rank-$k$ approximation of $\bA$ in terms of the spectral, trace and Frobenius error is given by $[\![\bA]\!]_k := \bU_k \bSigma_k \bV_k^\mathrm{T}$, where 
$\bSigma_k$ is a diagonal matrix of $k$ dominant singular values of $\bA$, and 
$\bU_k$ and $\bV_k$ contain the associated left and right singular vectors. In other words, we have
$$ [\![\bA]\!]_k = \arg \min_{\mathrm{rank}(\bB)=k} \| \bA - \bB  \|_{\xi},$$
where $\xi = 2, *$ or $\mathrm{F}$. Furthermore, $\| \bA - [\![\bA]\!]_k  \|_{\xi} = \sigma_{k+1}$ if $\xi = 2$, $\sum^m_{i=k+1} \sigma_i$ if $\xi = *$,  or $(\sum^m_{i=k+1} \sigma^2_i)^{\frac{1}{2}}$ if $\xi = \Frob$, where $\sigma_{1} \geq \sigma_2 \geq \hdots \geq \sigma_m$ denote the singular values of $\bA$. 

Obtaining the best rank-$k$ approximation can be computationally expensive and often becomes the bottleneck of the algorithm. In such case, one has to turn to alternative methods for computing a low-rank approximation, such as the randomized methods described below. 

\subsection{Randomized Singular Value Decomposition} \label{randSVDs}

A low-rank approximation of $\bA$ can be interpreted as reduction of the range of $\bA$ to a low-dimensional subspace $Q$ capturing the most of $\bA$'s action. 
The SVD approximation $[\![\bA]\!]_k$ corresponds to taking $Q$ as $\mathrm{range}(\bU_k)$.  A more efficient way is to take $Q = \mathrm{range}(\bA \bOmega^\mathrm{T})$, where $\bOmega\in \mathbb{R}^{l \times n}$ is an OSE~\cite{halko2011finding,woodruff2014sketching}. In this case, the optimal approximation is
\begin{equation} \label{eq:RSVD} 
[\![\bA]\!]^{\scriptscriptstyle\mathrm{(RSVD)}} := \arg \min_{\mathrm{range}(\bB) \subseteq Q}  \| \bA - \bB  \|_{\xi}.
\end{equation}
Then a rank-$k$ approximation of $\bA$ can be obtained by a truncated SVD of $[\![\bA]\!]^{\mathrm{\scriptscriptstyle(RSVD)}}$, which leads to the approximation $[\![\bA]\!]^{\scriptscriptstyle\mathrm{(RSVD)}}_k := [\![ [\![\bA]\!]^{\mathrm{\scriptscriptstyle(RSVD)}}]\!]_k$.
Notice that  $\bQ^\mathrm{T}[\![\bA]\!]^{\scriptscriptstyle\mathrm{(RSVD)}}_k = [\![\bQ^\mathrm{T}\bA]\!]_k$, where $\bQ$ is an orthonormal basis for $Q$. This observation constitutes the core for the RSVD algorithm (see \cref{alg:RSVD}) for the computation of $[\![\bA]\!]^{\scriptscriptstyle \mathrm{(RSVD)}}_k$. { \Cref{alg:RSVD} is same as \cite[Algorithm 8]{martinsson2020randomized}, with specifying the way of computing the SVD of a tall and skinny matrix $\bZ^\mathrm{T}$ by a QR factorization and the SVD of the R factor (see steps 4 and 5). The computation of $\bY$ can be effectively done with the procedure from~\cref{blockSRHT} using the block structure of $\bOmega$. 
The QR factorizations in steps 2 and 4 should be performed with TSQR~\cite{demmel2012communication} or other methods having low communication cost. 
The computational cost of~\cref{alg:RSVD} is dominated by computing $\bQ^\mathrm{T} \bA$ in step 3 and possibly $\bA \bOmega^\mathrm{T}$ in step 1. }

\small 
\begin{algorithm} 
	\caption{RSVD, based on~\cite[Algorithm 8]{martinsson2020randomized}} \label{alg:RSVD}
	\begin{algorithmic}[1] 		\Require{$m \times n$ matrix  $\bA$, $l \times n$ matrix $\bOmega$ with $l \ll n$, the target rank $k$.}
		\State{Compute $\bY = \bA \bOmega^\mathrm{T}$.}
		\State{Orthogonalize $\bY$ with QR factorization, and get $\bQ$.}
		\State{Compute matrix $\bZ=\bQ^\mathrm{T} \bA$.} 
		\State{Obtain a QR factorization $\bP \bR$ of $\bZ^\mathrm{T}$} 
		\State{Use SVD to compute the best rank-$k$ approximation $\tilde{\bU}_k \tilde{\bSigma}_k \tilde{\bV}_k^\mathrm{T} $ of  $\bR^\mathrm{T}$.}
		\State{Output factorization $[\![\bA]\!]^{\scriptscriptstyle\mathrm{(RSVD)}}_k = (\bQ \tilde{\bU}_k) \tilde{\bSigma}_k (\bP \tilde{\bV}_k)^\mathrm{T}$.}
	\end{algorithmic}
\end{algorithm}
\normalsize

Let us now characterize the accuracy of RSVD approximation. {It can be measured for instance by the quasi-optimality constant $\frac{\| \bA - [\![\bA]\!]^{\scriptscriptstyle\mathrm{(RSVD)}}_k \|_{\xi}}{\| \bA - [\![\bA]\!]_k  \|_{\xi}}-1$.} It is first shown that the accuracy of $[\![\bA]\!]^{\scriptscriptstyle\mathrm{(RSVD)}}_k$ is guaranteed if $Q = \mathrm{range}(\bA \bOmega^{\mathrm{T}})$ captures well the range of $\bA$, i.e., for some $d \geq k$ and $\varepsilon^* \leq \frac{1}{2}$ we have
\begin{equation} \label{eq:PiQA12}
	\| \bA -[\![\bA]\!]^{\scriptscriptstyle\mathrm{(RSVD)}}\|^2_\Frob \leq (1+\varepsilon^*) \| \bA - [\![\bA]\!]_d\|^2_\Frob. 
\end{equation}
Then, by the triangle inequality, 
\begin{equation} \label{eq:tineq}
\| \bA - [\![\bA]\!]^{\scriptscriptstyle\mathrm{(RSVD)}}_k \|_{\xi} \leq \| \bA - [\![\bA]\!]_k  \|_{\xi} + \| \bA -  [\![\bA]\!]^{\scriptscriptstyle\mathrm{(RSVD)}}\|_{\xi},
\end{equation}
where $\xi = 2, *$ or $\Frob$, we obtain 
\begin{equation} \label{eq:Nystrk}
\| \bA - [\![\bA]\!]^{\scriptscriptstyle\mathrm{(RSVD)}}_k \|_{\xi} \leq \| \bA - [\![\bA]\!]_k  \|_{\xi} + 2.5\| \bA -  [\![\bA]\!]_d\|_{\Frob}.
\end{equation}		 
This result guarantees, under the condition~\cref{eq:PiQA12}, the quasi-optimality of $[\![\bA]\!]^{\scriptscriptstyle\mathrm{(RSVD)}}_k $ with respect to the Frobenius norm. Furthermore, if $\bA$ has a fast enough singular value decay and $d$ is large enough,  so that the tail after $d$-th singular value of $\bA$, i.e. $(\sum^m_{i=d+1 } \sigma^2_i)^{\textstyle \frac{1}{2}} = \| \bA - [\![\bA]\!]_d\|_{\Frob}$ is small compared to the $k+1$-th singular value $\sigma_{k+1} = \| \bA - [\![\bA]\!]_k \|_{2}$ then $[\![\bA]\!]^{\scriptscriptstyle\mathrm{(RSVD)}}_k$ is almost as accurate as $[\![\bA]\!]_k$ with respect to all three norms. Moreover, increasing $d$ can make the quasi-optimality constant arbitrary close to zero.

It remains to obtain the conditions on $\bOmega$ such that~\cref{eq:PiQA12} holds with high probability. This can be done for instance with the results from~\cite{woodruff2014sketching}. Take $\varepsilon = \frac{4}{9} \varepsilon^*$. It follows from~[Lemma 45]\cite{woodruff2014sketching} and its proof that~\cref{eq:PiQA12} holds with probability at least $1-\delta$ if $\bOmega$ is an $(\frac{1}{3},\delta,d)$ OSE, and 
$$ \|\bV_d^\mathrm{T} \bOmega^\mathrm{T} \bOmega (\bA - [\![\bA]\!]_d)^\mathrm{T} \|^2_\Frob \leq \varepsilon \| \bA - [\![\bA]\!]_d \|^2_\Frob. $$
In turn the latter condition is satisfied with probability at least $1-\delta$ if $\bOmega$ is an $(\sqrt{\frac{\varepsilon}{d}}, \frac{\delta}{N},1)$ OSE, where $N = 2md+m+d$, as shown below. The OSE property of $\bOmega$ and the union bound argument guarantee that for given $N$ fixed vectors $\bz_i$, we have
		\begin{equation} \label{eq:vecnrm}
		(1-{\textstyle\sqrt{\frac{\varepsilon}{d}}})\| \bz_i \|^2_2 \leq \| \bOmega \bz_i  \|^2_2 \leq (1+{\textstyle\sqrt{\frac{\varepsilon}{d}}})\| \bz_i \|^2_2,\text{ for}~1\leq i \leq N
		\end{equation}
		with probability at least $1-\delta$. Take set $\{\bz_i\}$ composed of the columns of $\bV_d$ denoted by $\bx_i$, the columns of $(\bA - [\![\bA]\!]_d)^\mathrm{T}$ denoted by $\by_i$, and all the pairs $\bx_i+\by_j$ and $\bx_i-\by_j$. Then the relation~\cref{eq:vecnrm}, the parallelogram identity and the fact that $\bx_i^\mathrm{T} \by_j =0$, imply that
		$$ | \bx_i^\mathrm{T} \bOmega^\mathrm{T} \bOmega \by_j | \leq {\textstyle\sqrt{\frac{\varepsilon}{d}}}\|\bx_i \|_2 \|\by_j\|_2, \text{ for}~ 1\leq i \leq d,~ 1\leq j\leq m. $$
		Consequently, we have
		\small
		\begin{equation}
		\begin{split}
		\|\bV_d^\mathrm{T} \bOmega^\mathrm{T} \bOmega (\bA - [\![\bA]\!]_d)^\mathrm{T} \|^2_\Frob = \sum^d_{i=1} \sum^m_{j=1} |\bx_i^\mathrm{T} \bOmega^\mathrm{T} \bOmega \by_j|^2 
		\leq {\textstyle\frac{\varepsilon}{d}} \sum^d_{i=1} \sum^m_{j=1} \|\bx_i \|^2_2 \|\by_j\|^2_2 =  {\textstyle\frac{\varepsilon}{d}} \| \bV_d \|^2_\Frob \| \bA - [\![\bA]\!]_d \|^2_\Frob 
		\end{split}
		\end{equation}
		\normalsize
		with probability at least $1-\delta$. The proof is finished by noting that $\| \bV_d \|^2_\Frob = d$.  
		
		It is concluded that~\cref{eq:PiQA12} and as a consequence~\cref{eq:Nystrk} are satisfied with probability at least $1-2\delta$ if $\bOmega$ is an $(\frac{1}{3},\delta,d)$ OSE and $(\sqrt{\frac{\varepsilon}{d}}, \frac{\delta}{N},1)$ OSE. In turn, according to~\cref{thm:main} and~\cref{eq:SRHTbound}, this condition is satisfied by the block as well as the standard SRHT with $l = \mathcal{O}(d \log\!\frac{n}{\delta}\log\!\frac{m}{\delta})$ rows (taking $\varepsilon^*=\frac{1}{2}$,  $\varepsilon=\frac{2}{9}$). Whereas for Gaussian matrices the required number of rows to satisfy the aforementioned OSEs properties is somewhat lower: $l = \mathcal{O}(d \log{\frac{m}{\delta}})$. Although it has to be said that SRHT matrices in practice give similar results as Gaussian matrices~\cite{halko2011finding}. As can be seen from our experiments, this should also be the case for block SRHT.
		Moreover, we note that the condition $l = \mathcal{O}(d \log{\frac{n}{\delta}})$ for Gaussian matrices is still pessimistic. This overestimation is an artifact due to the use of a general analysis based solely on the OSE property.
		In reality, a Gaussian $\bOmega$ should satisfy~\cref{eq:PiQA12} with high probability if it has size $l = \mathcal{O}(d)$ with a small constant (say $2$ or $4$)~\cite{halko2011finding,tropp2017fixed}.

		\subsection{Nystr\"{o}m approximation}
		
		Although being more efficient than the deterministic SVD, the RSVD still can be computationally heavy, especially on distributed architectures, as it requires two passes over $\bA$, and a multiplication of $\bA$ by $\bQ$, which is a large dense matrix. Next we discuss improved algorithm that can circumvent these drawbacks. Assume that $\bA$ is positive semi-definite matrix.
		
		Notice that it can be computationally beneficial to change the norm $\| \cdot \|_\xi$ in~\cref{eq:RSVD} to its sketched estimate  $\| \bOmega \cdot \bOmega^\mathrm{T}\|_\xi$. The accuracy of such an estimation can be guaranteed thanks to the fact that $\bOmega$ is an OSE. 
		This leads to Nystr\"{o}m approximation $[\![\bA]\!]^{\mathrm{\scriptscriptstyle(Nyst)}}$ given below
		\begin{equation} \label{eq:Nystrom} 
		[\![\bA]\!]^{\mathrm{\scriptscriptstyle(Nyst)}} := \arg \min_{\substack{\mathrm{range}(\bB) \subseteq Q}}  \| \bOmega (\bA - \bB) \bOmega^\mathrm{T}  \|_{\xi},
		\end{equation}
		or in a more usual form~\cite{alaoui2015fast,tropp2017fixed,gittens2013revisiting,chiu2013sublinear,drineas2005nystrom}:
		$$ [\![\bA]\!]^{\mathrm{\scriptscriptstyle(Nyst)}} = (\bOmega \bA)^\mathrm{T}  (\bOmega \bA \bOmega^\mathrm{T})^{\dagger}  (\bOmega\bA), $$
		where $(\bOmega \bA \bOmega^\mathrm{T})^{\dagger}$ denotes the pseudo-inverse of $\bOmega \bA \bOmega^\mathrm{T}$.
		Then a rank-$k$ approximation of $\bA$ can be obtained by an SVD of $[\![\bA]\!]^{\mathrm{\scriptscriptstyle(Nyst)}}$, which leads to the approximation $[\![\bA]\!]^{\scriptscriptstyle\mathrm{(Nyst)}}_k := [\![ [\![\bA]\!]^{\mathrm{\scriptscriptstyle(Nyst)}}]\!]_k$. This way of obtaining a rank-$k$ approximation from $[\![\bA]\!]^{\scriptscriptstyle\mathrm{(Nyst)}}$ is referred to as the modified fixed-rank Nystr\"{o}m via QR~\cite{pourkamali2019improved,pourkamali2018randomized,tropp2017fixed}	
		{\Cref{alg:Nystrom} describes a way for computing $[\![\bA]\!]^{\mathrm{\scriptscriptstyle(Nyst)}}_k$ suited for distributed computing under 2D partitioning of $\bA$. The matrices $\bY$ and $\bOmega\bY$ can be computed with the procedure from~\cref{blockSRHT} using the block structure of $\bOmega$. 
			The QR factorization $\bZ = \tilde{\bQ} \bR$ in step 4 can be computed with TSQR.  Note that in step 6, instead of computing $\hat{\bU}_k$  as $(\bY \tilde{\bV}_k) \tilde{\bSigma}_k^{-1}$ we  could use  $\hat{\bU}_k = \tilde{\bQ} \tilde{\bU}_k$, which would provide more numerical stability but entail a larger computational cost. }{ \Cref{alg:Nystrom} needs only one pass over the matrix $\bA$, and does not involve any high-dimensional operations on $\bA$ except the computation of the sketch $\bY=\bA \bOmega^\mathrm{T}$, which implies its superiority over the standard SVD as well as randomized SVD~\cite{halko2011finding,tropp2017fixed}. In fact, the dominant computational cost of~\cref{alg:Nystrom} is associated with computing $\bY$ and $\bOmega\bY$ in steps 1 and 2, when $r$ is sufficiently large. } 
		\begin{algorithm} 
			\caption{Randomized Nystr\"{o}m approximation} \label{alg:Nystrom}
			\label{algo_proj_2d}
			\begin{algorithmic}[1] 		\Require{$n \times n$ matrix  $\bA$, $l \times n$ matrix $\bOmega$ with $l \ll n$, the target rank $k$.}
				\State{Compute $\bY = \bA \bOmega^\mathrm{T}$.}
				\State{Obtain a Cholesky factor $\bC$ of $\bOmega \bY$. }
				\State{Compute $\bZ = \bY \bC^{-1}$ with backward substitution.}
				\State{Obtain the R factor $\bR$ of $\bZ$ (with TSQR or similar algorithm).}	
				\State{Use SVD to compute the best rank-$k$ approximation $\tilde{\bU}_k \tilde{\bSigma}_k \tilde{\bV}_k^\mathrm{T} $ of  $\bR$.}
				\State{Compute $\hat{\bU}_k = (\bY \tilde{\bV}_k) \tilde{\bSigma}_k^{-1}$.} 	
				\State{Output factorization $[\![\bA]\!]^{\mathrm{\scriptscriptstyle(Nyst)}}_k = \hat{\bU}_k \tilde{\bSigma}^2_k  \hat{\bU}_k^\mathrm{T}$.}
			\end{algorithmic}
		\end{algorithm}
		\begin{remark} \label{rmk:cholesky}
			\color{black}
			The matrix $\bOmega \bA \bOmega^\mathrm{T}$ can be rank-deficient, for instance, if $\bA$ or $\bOmega$ have lower rank than $l$, which will cause a problem for obtaining a Cholesky factorization in step 2.   In this case, a remedy can be to compute an SVD instead of the Cholesky factorization, and take $\bC$ as a square root of $\bOmega \bY$ in SVD form, that then can be used for the pseudo-inversion in step 3. Another possibility is to make $\bA$ full-rank by using shifting as in~\cite{li2017algorithm}. 
			\color{black}
		\end{remark}	
		
		Let us now characterize the accuracy of Nystr\"{o}m approximation. Notice the following identity~\cite{gittens2011spectral}:
		$$ \bA - [\![\bA]\!]^{\mathrm{\scriptscriptstyle(Nyst)}}= (\bA^{\frac{1}{2}} - [\![\bA^{\frac{1}{2}}]\!]^{\mathrm{\scriptscriptstyle(RSVD)}})^\mathrm{T} (\bA^{\frac{1}{2}} - [\![\bA^{\frac{1}{2}}]\!]^{\mathrm{\scriptscriptstyle(RSVD)}}).  $$
		By combining this observation with the derived earlier results on RSVD with $\bA \leftarrow \bA^{\frac{1}{2}}$, we obtain that 
		\begin{equation} \label{eq:Nyst}
		\begin{split}
		&\|\bA - [\![\bA]\!]^{\mathrm{\scriptscriptstyle(Nyst)}} \|_{*} = \|\bA^{\frac{1}{2}} - [\![\bA^{\frac{1}{2}}]\!]^{\mathrm{\scriptscriptstyle(RSVD)}} \|^2_\Frob  \leq (1+\varepsilon^*)\|\bA^{\frac{1}{2}} - [\![\bA^{\frac{1}{2}}]\!]_d \|^2_\Frob  = (1+\varepsilon^*) \|\bA - [\![\bA]\!]_d \|_{*}
		\end{split}
		\end{equation}	
		holds with probability at least $1-2\delta$ if $\bOmega$ is $(\frac{1}{3},\delta,d)$ OSE and $(\sqrt{\frac{\varepsilon}{d}}, \frac{\delta}{N},1)$ OSE.
		It is then noticed that~\cref{eq:Nyst} also implies the accuracy of the truncated approximation $[\![\bA]\!]^{\scriptscriptstyle\mathrm{(Nyst)}}_k$ due to the following consequence of the triangle inequality (see for instance~\cite[Proposition A.6]{tropp2017practical}):
		\begin{equation} \label{eq:tineq}
		\| \bA - [\![\bA]\!]^{\scriptscriptstyle\mathrm{(Nyst)}}_k \|_{\xi} \leq \| \bA - [\![\bA]\!]_k  \|_{\xi} + 2\| \bA -  [\![\bA]\!]^{\scriptscriptstyle\mathrm{(Nyst)}}\|_{\xi},
		\end{equation}
		where $\xi = 2$ or $*$, so that we have by~\cref{eq:Nyst},
		\begin{equation} \label{eq:Nystrk}
		\| \bA - [\![\bA]\!]^{\scriptscriptstyle\mathrm{(Nyst)}}_k \|_{\xi} \leq \| \bA - [\![\bA]\!]_k  \|_{\xi} + 3\| \bA -  [\![\bA]\!]_d\|_{*}.
		\end{equation}
		Similarly to RSVD approximation, this relation guarantees the quasi-optimality of $[\![\bA]\!]^{\scriptscriptstyle\mathrm{(Nyst)}}_k $ with respect to the trace norm. Furthermore,  $[\![\bA]\!]^{\scriptscriptstyle\mathrm{(Nyst)}}_k$ is almost as accurate as $[\![\bA]\!]_k$ with respect to both the trace norm and the spectral norm, if the tail after $d$-th singular value of $\bA$, i.e. $\sum^n_{i=d+1 } \sigma_i$ is small compared to the $k+1$-th singular value $\sigma_{k+1}$.
		
		\subsection{Single-view approximation of non-psd matrix}
The Nystr\"{o}m method is applicable only when $\bA$ is positive semi-definite.  Next, we describe a single-view algorithm that works with general matrices. In recent years several such algorithms have been proposed~\cite{halko2011finding,upadhyay2016fast,tropp2017practical,tropp2019streaming,kannan2014principal}. The single-view approximation from~\cite{tropp2019streaming} involves four sketching matrices. Two of them, $\bOmega \in \mathbb{R}^{l \times n}$, and $ \bGamma \in \mathbb{R}^{l \times m}$, are used to construct the approximation subspaces $Q = \mathrm{range}(\bA \bOmega^\mathrm{T})$ and $P = \mathrm{range}(\bA^\mathrm{T} \bGamma^\mathrm{T})$ capturing the actions of $\bA$ and $\bA^\mathrm{T}$. 
Whereas the other two, $\bPhi \in \mathbb{R}^{s \times m}$, and $\bPsi \in \mathbb{R}^{s \times n}$, with $l \leq s$, provide the ``core sketch'' of $\bA$. In our {projection-based} interpretation, the ``core sketch'' corresponds to estimation of the norm $\| \cdot \|_\xi$ of the residual  by $\| \bPhi \cdot \bPsi^\mathrm{T} \|_\xi$. Thus, we arrive to the following low-rank approximation of $\bA$:
\begin{equation} \label{eq:spass}
[\![\bA]\!]^{\mathrm{\scriptscriptstyle(sRSVD)}} = \arg \min_{\substack{\mathrm{range}(\bB) \subseteq Q, \\
		\mathrm{range}(\bB^\mathrm{T}) \subseteq P} } \| \bPhi (\bA - \bB) \bPsi^\mathrm{T}  \|_{\xi}.
\end{equation}
{Notice that $[\![\bA]\!]^{\mathrm{\scriptscriptstyle(sRSVD)}}$ is given as $\bQ (\bPhi \bQ)^\dagger (\bPhi \bA \bPsi^\mathrm{T}) (\bPsi \bP)^\dagger \bP^\mathrm{T}$, where $\bQ$ denotes an orthogonal basis for $Q$ and $\bP$ denotes an orthogonal basis for $P$. 
Similarly as in the case of RSVD and Nystr\"{o}m approximations, the rank-$k$ approximation of $\bA$ can be obtained by truncating $[\![\bA]\!]^{\mathrm{\scriptscriptstyle(sRSVD)}}$ with SVD, which provides $[\![\bA]\!]^{\mathrm{\scriptscriptstyle(sRSVD)}}_k= [\![[\![\bA]\!]^{\scriptscriptstyle\mathrm{(sRSVD)}}]\!]_k$~\cite{tropp2019streaming}.} Notice that $[\![\bA]\!]^{\mathrm{\scriptscriptstyle(sRSVD)}}_k$ is now given as $\bQ [\![\bC]\!]_k \bP^\mathrm{T}$, where $\bC = (\bPhi \bQ)^\dagger (\bPhi \bA \bPsi^\mathrm{T}) (\bPsi \bP)^\dagger$. This observation leads to~\cref{alg:sRSVD} for the computation of $[\![\bA]\!]^{\scriptscriptstyle\mathrm{\scriptscriptstyle(RSVD)}}_k$. Again, in step 2 one can use TSQR algorithm for the efficient QR factorization on distributed architectures. 
{As in Nystr\"{o}m method, the computational cost of~\cref{alg:sRSVD} is dominated by the applications of sketching matrices.}

\begin{algorithm} 
	\caption{Single-view RSVD~\cite[Algorithm 17]{martinsson2020randomized}} \label{alg:sRSVD}
	\begin{algorithmic}[1] 		\Require{$m \times n$ matrix  $\bA$, matrices $\bOmega, \bGamma$ with $l$ rows,  and $\bPhi, \bPsi$ with $s$ rows, with $l \leq s \ll m$, the target rank $k$.}
		\State{Compute $\bX =\bA^\mathrm{T} \bGamma^\mathrm{T}$, $\bY = \bA \bOmega^\mathrm{T}$ and $\bZ = \bPhi \bA \bPsi^\mathrm{T}$.}
		\State{Orthogonalize $\bX$ and $\bY$ with a QR factorization and obtain $\bQ$ and $\bP$.}
		\State{Compute $\bPhi \bQ$ and $\bPsi \bP$.}
		\State{Compute the core matrix $\bC = (\bPhi \bQ)^\dagger \bZ (\bPsi \bP)^\dagger$ with least-squares solves.}
		\State{Use SVD to compute the best rank-$k$ approximation $\tilde{\bU}_k \tilde{\bSigma}_k \tilde{\bV}_k^\mathrm{T} $ of  $\bC$.}
		\State{Output factorization $[\![\bA]\!]^{\mathrm{\scriptscriptstyle(sRSVD)}}_k = (\bQ \tilde{\bU}_k) \tilde{\bSigma}_k (\bP \tilde{\bV}_k)$.}
	\end{algorithmic}
\end{algorithm}

To characterize the accuracy of  $[\![\bA]\!]^{\mathrm{\scriptscriptstyle(sRSVD)}}_k$ we shall assume that $\bPhi, \bPsi$ satisfy both  $(\varepsilon, \delta,l)$ OSE and $(\varepsilon, \frac{\delta}{n},1)$ OSE properties with $\varepsilon \leq \frac{2}{9}$. Then these sketching matrices are $\varepsilon$-embeddings for $Q$ and $P$, i.e. they satisfy~\cref{eq:isometry} taking $\bOmega \leftarrow \bPhi, V \leftarrow Q$ or $\bOmega \leftarrow \bPsi,V \leftarrow P$, simultaneously, with probability at least $1-2\delta$. Moreover, let  $\bPi_{Q}$ and $\bPi_{P}$ denote the orthogonal projectors onto $Q$ and $P$.  Then $\bPsi$ satisfies the $\varepsilon$-embedding property for every subspace spanned by a row of $\bA - \bPi_{Q} \bA \bPi_{P}$ simultaneously with probability at least $1-\delta$, and $\bPhi$ satisfies the $\varepsilon$-embedding property for every subspace spanned by a column of $(\bA - \bPi_{Q} \bA \bPi_{P})\bPsi^\mathrm{T}$ simultaneously with probability at least $1-\delta$, so that 
$$\| \bPhi(\bA - \bPi_{Q} \bA \bPi_{P})\bPsi^\mathrm{T}\|_\Frob \leq \sqrt{1+\varepsilon}\| (\bA - \bPi_{Q} \bA \bPi_{P})\bPsi^\mathrm{T}\|_\Frob \leq (1+\varepsilon)\| \bA - \bPi_{Q} \bA \bPi_{P}\|_\Frob $$
holds with probability at least $1-2\delta$.
Thus, we have with probability at least $1-4\delta$,
	\begin{align*}
	&\| \bA - [\![\bA]\!]^{\mathrm{\scriptscriptstyle(sRSVD)}}  \|_\Frob \leq \| \bPi_{Q}\bA \bPi_{P} - [\![\bA]\!]^{\mathrm{\scriptscriptstyle(sRSVD)}} \|_\Frob  +  \| \bA - \bPi_{Q} \bA \bPi_{P}\|_\Frob \\
	&  \leq \frac{1}{{1-\varepsilon}} \| \bPhi (\bPi_{Q}\bA \bPi_{P} - [\![\bA]\!]^{\mathrm{\scriptscriptstyle(sRSVD)}})\bPsi^\mathrm{T}\|_\Frob  +  \| \bA - \bPi_{Q} \bA \bPi_{P}\|_\Frob \\
	&  \leq \frac{1}{{1-\varepsilon}} \| \bPhi (\bA - [\![\bA]\!]^{\mathrm{\scriptscriptstyle(sRSVD)}})\bPsi^\mathrm{T}\|_\Frob  +  2.6\| \bA - \bPi_{Q} \bA \bPi_{P}\|_\Frob \\
	&  \leq \frac{1}{{1-\varepsilon}} \| \bPhi (\bA - \bPi_{Q}\bA \bPi_{P})\bPsi^\mathrm{T}\|_\Frob  +  2.6\| \bA - \bPi_{Q} \bA \bPi_{P}\|_\Frob \\
	&  \leq  4.2\| \bA - \bPi_{Q} \bA \bPi_{P}\|_\Frob \leq  4.2(\| \bA - \bPi_{Q} \bA \|_\Frob + \| \bA -  \bA\bPi_{P} \|_\Frob),
	\end{align*}
which in turn guarantees the accuracy of $[\![\bA]\!]^{\mathrm{\scriptscriptstyle(sRSVD)}}$ if $Q$ and $P$ capture well the actions of $\bA$ and $\bA^\mathrm{T}$.

Clearly, $\bPi_{Q} \bA  =  [\![\bA]\!]^{\scriptscriptstyle\mathrm{(RSVD)}}$ and $\bA \bPi_{P} = ([\![\bA^\mathrm{T}]\!]^{\scriptscriptstyle\mathrm{(RSVD)}})^\mathrm{T}$, where $[\![\bA^\mathrm{T}]\!]^{\scriptscriptstyle\mathrm{(RSVD)}}$ is the RSVD approximation of $\bA^\mathrm{T}$, associated with the sketching matrix $\bGamma$. This suggests that the single-view approximation~\cref{eq:spass} should have a similar quality as the RSVD approximation~\cref{eq:RSVD}. By combining the above consideration  with the results on RSVD from~\cref{randSVDs} and the triangle inequality  we obtain that 
\small 
	\begin{align*} 
	\| \bA - [\![\bA]\!]^{\mathrm{\scriptscriptstyle(sRSVD)}}_k \|_{\xi}  &\leq  \| \bA - [\![\bA]\!]^{\mathrm{\scriptscriptstyle(sRSVD)}}_k \|_{\xi}+	2\| \bA - [\![\bA]\!]^{\mathrm{\scriptscriptstyle(sRSVD)}} \|_{\Frob} \\&\leq \| \bA - [\![\bA]\!]^{\mathrm{\scriptscriptstyle(sRSVD)}}_k \|_{\xi}+16.8\sqrt{1+2.25\varepsilon} \|\bA - [\![\bA]\!]_d \|_\Frob
	\end{align*}
	\normalsize
holds with probability at least $1-6\delta$ for $\xi = 2,*$ or $\Frob$, given that $\bPhi, \bPsi$ satisfy both $(\varepsilon, \delta,l)$ OSE and $(\varepsilon, \frac{\delta}{n},1)$ OSE properties, and $\bOmega,\bGamma$ satisfy both $(\frac{1}{3},\delta,d)$ OSE and $(\sqrt{\frac{\varepsilon}{d}}, \frac{\delta}{N},1)$ OSE properties, with $\varepsilon \leq \frac{2}{9}$ and $N=2nd+n+d$. According to~\cref{thm:main}, these conditions can be satisfied by  $\bOmega,\bGamma$ that are block SRHT matrices of size $l = \mathcal{O}(d \log^2\!\frac{n}{\delta})$, and $\bPhi, \bPsi$ that are block SRHT matrices of size $s = \mathcal{O}(d \log^3\!\frac{n}{\delta})$. 
		\section{Numerical experiments} \label{experiments}
		For numerical experiments, we chose the Nystr\"{o}m approximation as a representative application. 
		The validation of block SRHT is done through comparison with Gaussian embeddings. In the plots, BSRHT refers to block SRHT. The comparison with standard SRHT is impertinent since SRHT matrices are not that well scalable as Gaussian matrices and have no better accuracy~\cite{yang2015implementing}.

		\subsection{Nystr\"{o}m approximation}	
		This experiment was executed with Julia programming language version 1.7.2 along with the Distributed.jl
		and DistributedArrays.jl packages for parallelism. We used 2
		nodes Intel Skylake 2.7GHz (AVX512) having 48 available cores and 180 MB of RAM each. In this experiment we used only 32 cores on each node.
		{As input data we took the MNIST or YearPredictionMSD datasets~\cite{726791,Bertin-Mahieux2011}.}
		The radial basis function $e^{-\|x_i - x_j\|^2 / \sigma^2}$ was used to build a dense positive definite matrix $\bA$ of size $n \times n$ from $n$ rows of the input data. The parameter
		$\sigma$ was chosen as $100$ for the MNIST dataset {\color{black}and $10^4$ as well as $10^5$ for the YearPredictionMSD dataset.}
		The dimension $n$ was taken as $65536$. The matrix $\bA$ has been uniformly distributed
		on a square grid of $8 \times 8$ processors. 
		In all the experiments, the local matrices $\bOmega^{(i)}$ on each processor were generated with a
		seeded random number generator with a low communication cost. 
		\begin{figure}[htpb]
			\begin{subfigure}[b]{0.2\textwidth}
				\input{accuracy_pot16_srht.tex}
				\subcaption{Dataset \texttt{mnist}}
			\end{subfigure}
			\begin{subfigure}[b]{0.3\textwidth}
				\input{accuracy_pot16_year_1e4.tex}
				\subcaption{Dataset \texttt{year} with $\sigma = 10^4$}
			\end{subfigure}
			\begin{subfigure}[b]{0.3\textwidth}
				\input{accuracy_pot16_year_1e5.tex}
				\subcaption{Dataset \texttt{year} with $\sigma = 10^5$}
			\end{subfigure}
			
			\caption{Trace error $\|\bA-[\![\bA]\!]^{\mathrm{\scriptscriptstyle(Nyst)}}_k\|_{{*}} / \|\bA\|_{{*}}$ using BSRHT.}%
			\label{fig:accuracy}
		\end{figure}
		
		\Cref{fig:accuracy} depicts the convergence of the error of the low-rank approximation obtained with 
		\cref{algo_proj_2d} taking $\bOmega$ as a block SRHT. {\color{black}The results for Gaussian $\bOmega$ are practically identical and therefore are not displayed.}
		
		In this numerical experiment, the error is measured with the trace norm.  Different
		sketching sizes $l$ were tested. For each pair of parameters $(l,k)$ $20$ different
		approximations were computed for each type of $\bOmega$,  in order to have the $95\%$
		confidence interval. Nevertheless this interval is not displayed as it is too small to
		be visible. 
		\Cref{fig:timings} gives runtime characterization. In particular we depict the runtime spent
		on computing $\bY = \bA \bOmega^\mathrm{T}$ and $\bOmega \bY$ in steps $1$ and {$2$} of \cref{algo_proj_2d}.
		These operations will dominate the overall computational cost, when the block size is large enough.
		Nevertheless the reader should be aware that TSQR and the SVD of $\bR$ (step 4 and 5) are also important,
		especially when the sampling size is close to the block size.
		The parameter $k$ is not involved in steps 1 and 3 hence not mentioned in \cref{fig:timings}.
		
		\begin{figure}[ht]
			\centering
			\input{runtime_julia.tex}
			\caption{Runtimes of computing $\bY = \bA \bOmega^\mathrm{T}$ and $\bOmega \bY$ in \cref{algo_proj_2d} for different sampling sizes.}
			\label{fig:timings}
		\end{figure}
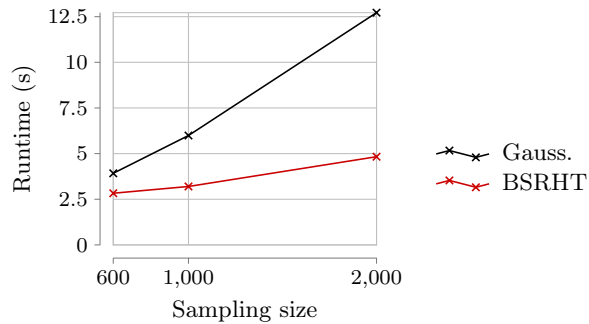
		
		\color{black}
		According to~\cref{fig:timings}, the runtime of the Gaussian sampling is up to 2.5 times higher and grows faster with $l$ than the runtime taken by block SRHT. Note that for block SRHT the local computation cost is independent of $l$, hence
		the slope comes only from the reductions in steps $1$ and $2$. On the other hand, the Gaussian sampling involves local computations with linear dependency in $l$, in addition to these reductions.
		\color{black}
		
		\subsection{Cost of application to tall-and-skinny matrix}	
		
		Next we investigate the performance of block SRHT on larger scale. For this we consider a product of $\bOmega$ with a tall-and-skinny matrix $\bV$, for instance in the context of solving an overdetermined least-squares problem. The same computing environment is used as in the previous experiment, involving now up to 32 nodes and using C99/MPI instead of Julia. The code was compiled using IntelMPI C compiler version 20.0.2 and sequential MKL 20.0.2 with option ILP64. The library FFTW3 used has Intel-specific routines. There is therefore up to 1536 cores available. In this way we generated a random matrix $\bV$ with $d=200$ columns and a variable number $n$ of rows. This matrix was distributed among a variable number $p$ of processors with block rowwise partitioning. Then $\bV$ was multiplied by either a Gaussian or block SRHT matrix $\bOmega$ with $l = 2000$ rows using~\cref{eq:blockSRHTappl}. In all experiments, the local $\bOmega^{(i)}$ matrices on each processor were generated with a seeded random number generator with negligible communication cost. \Cref{fig:StrongScala10E6} presents a strong scalability test for $n=10^7$. We see that the block SRHT provides an overall speedup by a factor of more than $2.5$ over the Gaussian matrices, while demonstrating as good scalability when $p \leq 384$. For larger $p$, however, the reduction operation starts to dominate, which reduces the gain in efficiency. We observe a variability in the MPI\_Allreduce operation on larger number of processors for both Gaussian and block SRHT algorithms. However the compute times for both algorithms scale well when increasing the number of processors up to $p = 1536$.   
		\begin{figure}[ht]
			\centering
			\input{B.tex}
			\caption{Strong scalability runtimes associated with computing $\bOmega \bV$ with $n = 10^7$ and $l = 2000$, versus $p$. ``Gauss. total'' and ``BSRHT total'' correspond to the overall runtimes, whereas ``Gauss. local'' and ``BSRHT local'' stand for the max per-processor runtimes taken by local multiplications.}
			\label{fig:StrongScala10E6}
		\end{figure}
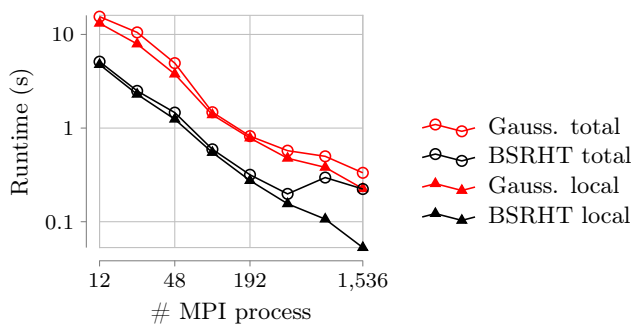
		\begin{figure}[ht]
			\centering
			\input{C.tex}
			\caption{Strong scalability runtimes associated with computing $\bOmega \bV$ with $n = 10^8$ and $l = 2000$, versus $p$. ``Gauss. total'' and ``BSRHT total'' correspond to the overall runtimes, whereas ``Gauss. local'' and ``BSRHT local'' stand for the max per-processor runtimes taken by local multiplications.}
			\label{fig:StrongScala100E6}
		\end{figure}
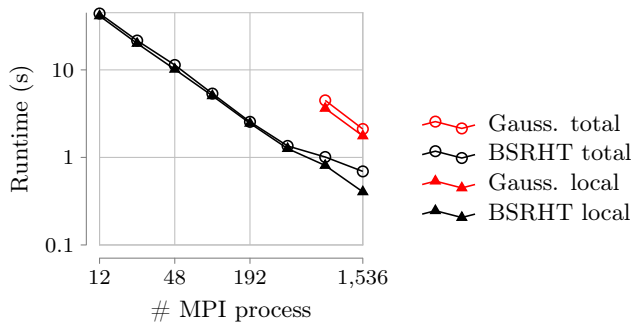
		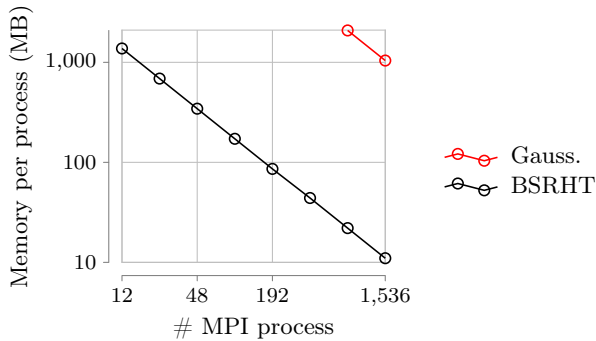
\begin{figure}[ht]
			\centering
			\input{100E6_Mem.tex}
			\caption{Max per-processor memory needed for computing  $\bOmega \bV$ with $n = 10^8$ and $l = 2000$, versus $p$.}
			\label{fig:MemStrongScala100E6}
		\end{figure}
		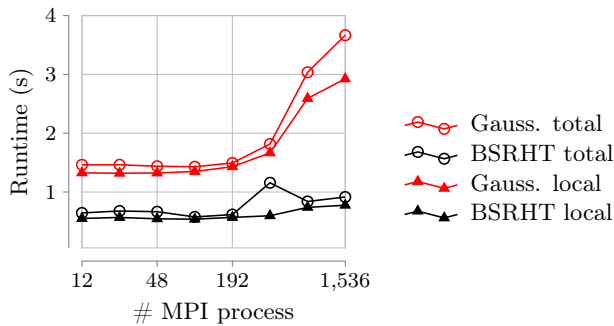
\begin{figure}[ht] 
			\centering
			\input{A.tex}
			\caption{Weak scalability runtimes associated with computing $\bOmega \bV$ with $n = 10^5 \times p $ and $l = 2000$, versus $p$. ``Gauss. total'' and ``BSRHT total'' correspond to the overall runtimes, whereas ``Gauss. local'' and ``BSRHT local'' stand for the max per-processor runtimes taken by local multiplications.}
			\label{fig:WeakScala}
		\end{figure}
		\Cref{fig:StrongScala100E6} shows a strong scalability test for a higher dimension $n=10^8$. Again we see a great scalability of block SRHT for $p \leq 384$. For Gaussian matrices, on the other hand, we revealed issues with reaching the memory limit needed to store $\bOmega^{(i)}$ which made its application on $p \leq 384$ processors infeasible. In principle, this problem can be overcome by generating $\bOmega^{(i)}$ blockwise and applying the blocks to $\bV^{(i)}$ "on the fly". This however entails a dramatic increase in runtime and therefore is omitted in comparison. On the other hand, for block SRHT we do not have any memory problems\footnote{To reduce the memory consumption, local matrices $\bV^{(i)}$ are multiplied by $\bOmega^{(i)}$ in blocks of $20$ columns.}. To quantify the advantage of block SRHT in such context, in~\Cref{fig:MemStrongScala100E6} we provide the memory consumption of the Gaussian and block SRHT matrix. We see that in this sense the reduction in computational cost is indeed drastic. Finally, \Cref{fig:WeakScala}~provides a weak scalability test using $n = rp$ where $r = 10^5$. We again see a reduction in runtime of about $2.5$ and good scalability for block SRHT up to using $p = 1536$ processors, similar as in the strong scalability test.

		\section{Proof of the main theorem} \label{mainproof}
		Before providing the proof for~\cref{thm:main}, let us first motivate the chosen proof path. 
		Le $\bV$ be a fixed $n \times d$ matrix with orthonormal columns, partitioned using block rowwise partitioning with $p$ blocks $\bV^{(i)}$ of size $r \times d$. The statement of the theorem can then be proven by showing that the singular values of $\bOmega \bV$ belong to the interval $[\sqrt{1- \varepsilon},\sqrt{1+\varepsilon}]$ with probability at least $1-\delta$.
		
		Assume for a moment, that $\bR$ in~\cref{eq:blockSRHT} is a uniform sampling matrix \textit{without replacement}. Notice that the random sampling of rows without replacement and then flipping their signs is equivalent to first flipping the signs and then sampling. By using this consideration, the expression~\cref{eq:blockSRHTappl} can be developed further as \small
		\begin{equation} \label{eq:blockSRHTappl2}
		\begin{split}
		\bOmega \bV &= \sqrt{\frac{r}{l}} \sum^p_{i=1} \left ( \widetilde{\bD}^{(i)} \bR \bH \bD^{(i)} \bV^{(i)}\right )
		= \sqrt{\frac{r}{l}} \sum^p_{i=1} \left ( \bR \bhD^{(i)} \bH \bD^{(i)} \bV^{(i)} \right ) \\ 
		&= \sqrt{\frac{r}{l}} \bR \sum^p_{i=1} \left (\bhD^{(i)} \bH \bD^{(i)} \bV^{(i)} \right ) = \sqrt{\frac{r}{l}} \bR\bW\bV,
		\end{split}
		\end{equation}
		\normalsize
		where $\bhD^{(i)}$ are $r \times r$ diagonal matrices with Rademacher random variables $\pm1$ on the diagonal, and   $\bW = \left [ \bhD^{(1)}\bH \bD^{(1)}, \bhD^{(2)} \bH \bD^{(2)}, \hdots, \bhD^{(l)} \bH \bD^{(l)} \right ]$. 
		
		Looking at~\cref{eq:blockSRHTappl}, one can detect many similarities of $\bOmega$ with standard SRHT matrix. Consequently, in order to argue that $\bOmega \bV$ is approximately orthonormal, the first thing to try should be to follow the steps from~\cite{tropp2011improved} in the analysis of the original SRHT. In this case the proof recipe would be as follows. First, it could be shown that the matrix $\bW$ with high probability homogenizes the rows  of $\bV$. This result then would allow the Matrix Chernoff concentration inequality from~\cite{tropp2011improved} to be applied to show that $\bW\bV$ and $\sqrt{\frac{r}{l}}\bR \bW \bV$ have approximately equal minimal and maximal singular values. With these results, it would remain to show that with high probability $\bW \bV$ is approximately orthonormal.  This, however, can be cumbersome or even impossible in some situations. Think, for example, of the situation when $r < d$.  Therefore, we will assume that $\bR$ is a uniform sampling matrix \textit{with replacement} and use the following trick.  For better presentation define parameters $\varepsilon^* = \frac{15}{16} \varepsilon $ and $\delta^* = \frac{\delta}{5}$.

		Recall that the sampling matrix $\bR$ restricts a vector $\bx = (x_1, \hdots, x_{r})$ to $l$ coordinates, i.e., we have 
		$$\bR \bx = (x_{i_1}, \hdots, x_{i_l}),~\text{with }1 \leq i_1, \hdots, i_l \leq r.$$ 
		The (multi-)set of indices $\{ i_1, \hdots, i_l \}$ is  a uniform random sample of $ \{1, \hdots, r \}$ with replacement. Notice that such sampling of indices is equivalent to the sampling uniformly at random with replacement from $\{1, \hdots, 1, 2, \hdots 2, \hdots, r, \hdots, r \}$ containing $K =  \lceil{ 10^4 \frac{n^2}{r\delta^*}}\rceil$ copies of each index. 
		This observation implies that the sampling matrix $\bR$ satisfies the identity
		$\bR \bH = \bhR [  \bH~\bH \hdots  \bH ]^\mathrm{T} =\bhR \bhH,$
		where $\bhR$ is  uniform sampling, with replacement, matrix of size $l \times rK$, and $\bhH$ is a block matrix with $K$ blocks of rows, each being equal to $\bH$. For a vector $\bx = (x_1, \hdots, x_{r K})$, matrix $\bhR$ satisfies
		\begin{equation} \label{eq:hRsampling}
		\bhR \bx = (x_{i_1}, \hdots, x_{i_l}),~\text{with }1 \leq i_1, \hdots, i_l \leq r K,
		\end{equation}
		where the indices $\{ i_1, \hdots, i_l \}$ are drawn uniformly at random \textit{with replacement} from $ \{1, \hdots, r K \}$. 
		Let $\mathcal{S}$ denote the event when $i_1, \hdots, i_l$ in~\cref{eq:hRsampling} are all disjoint indices. 
		\begin{lemma} \label{thm:Slemma}
			$\mathcal{S}$ occurs with probability at least $1-\delta^*$.
			\begin{proof}
				There are in total $\frac{(rK)^l}{l!}$ ways to select $l$ elements from a $(rK)$-element set and $\binom{rK}{l} = \frac{{rK(rK-1) \hdots (rK-l+1)}}{l!}$ ways to select $l$ disjoint elements. Consequently, we have 	
				\small
				\begin{align*}
				\mathbb{P}(S)= \prod^l_{i=1} (1-\frac{i-1}{rK} ) \geq  (1-\frac{l}{rK} )^l \geq 1-\frac{l^2}{rK} \geq 1- \delta^*.
				\end{align*}
				\normalsize		
			\end{proof}
		\end{lemma}  
		The goal will be to bound the singular values of $\bOmega \bV$ under the condition $\mathcal{S}$. The overall probability of success, then will follow by the union bound argument. Next is assumed that $\mathcal{S}$ is occurring. 
		Notice that, in this case, matrix $\bhR$ is equivalent to the matrix that samples the entries uniformly at random and \textit{without replacement}.
		Then, using the same arguments as in~\cref{eq:blockSRHTappl2}, we have the following expression for the product $\bOmega \bV$:
		\small
		\begin{equation} \label{eq:blockSRHTappl3}
		\begin{split}
		\bOmega \bV =  \sqrt{\frac{r}{l}} \sum^p_{i=1} \left ( \widetilde{\bD}^{(i)} \bhR \bhH \bD^{(i)} \bV^{(i)} \right ) 
		= \sqrt{\frac{r}{l}} \bhR \bhW \bV,
		\end{split}
		\end{equation}
		\normalsize
		where $\bhW$ is a block matrix composed of $K \times p$ blocks, with the $(j,i)$-th block being $\bhD^{(i,j)} \bH \bD^{(i)}$, where $\widehat{\bD}^{(i,j)}$ are diagonal matrices with entries i.i.d. Rademacher random variables $\pm1$. Unlike $\bW \bV$, the matrix $\bhW \bV$ (rescaled by $1/\sqrt{K}$) for sufficiently large $K$ can be  proven to be approximately orthonormal with high probability. We are ready to establish the proof of~\cref{thm:main}. 
		
		Notice that the condition in~\cref{thm:main} implies that
		\begin{equation} \label{eq:kprov}
		n \geq l \geq 3.2 {\varepsilon^*}^{-2} (\sqrt{d}+\sqrt{8 \log(r K/\delta^*)})^2 \log (d/\delta^*).
		\end{equation}
		
		In~\cref{thm:lemma1} is shown that, given $\mathcal{S}$, the matrix $\bhW \bV$ has rows with equilibrated norms.  
		
		\begin{proposition}\label{thm:lemma1}
			Given $\mathcal{S}$. The rows ${\bphi^{(j)}}$ of $\bhW \bV$ satisfy $$\mathbb{P} \left (  \max_{j=1, \hdots, r K}  \| \bphi^{(j)} \|_2 \leq {\textstyle \sqrt{\frac{d}{r}}+\sqrt{\frac{8\log(r K/\delta^*)}{r}} }\right ) \geq 1 - \delta^*.$$
			\begin{proof}
				Notice that, each row of $\bhW$ has entries that are i.i.d Rademacher random variables rescaled by $1/\sqrt{r}$.  Consequently, we have $$\bphi^{(j)} =  {\bxi^{(j)}}^\mathrm{T} \bV/\sqrt{r}, $$
				where $\bxi^{(j)}$ is a Rademacher vector.
				Define convex function $f(\bx) = \| \bx^\mathrm{T} \bV /\sqrt{r}\|_{{2}}$. Observe that $f(\bx)$ satisfies the Lipschitz bound:
				$$\forall \bx, \by,~|f(\bx) - f(\by)|\leq  \|\bx -\by \|_{{2}} \|\bV/{\sqrt{r}}\|_{{2}} =  \| \bx -\by\|_{{2}}/\sqrt{r}. $$
				This allows to apply the Rademacher tail bound 
				\begin{equation} \label{eq:lemma1}
				\mathbb{P} \left (f({\bxi^{(j)}}) \geq  \mathbb{E}f({\bxi^{(j)}})+t/\sqrt{r} \right) \leq \exp (-t^2/8),~\forall t \geq 0.
				\end{equation} 
				Observe that 
				$\mathbb{E}(f({\bxi^{(j)}})) \leq (\mathbb{E}(f({\bxi^{(j)}})^2))^{\frac{1}{2}} =\|\bV/{\sqrt{r}} \|_\mathrm{F} \leq \sqrt{d/{r}}$. The statement of the lemma follows by combining this relation with~\cref{eq:lemma1} with $t=\sqrt{8 \log{\frac{1}{\delta}}}$ and using the union bound argument.
			\end{proof}
		\end{proposition} 	
		
		In~\cref{thm:singWV} is proven that  $\frac{1}{\sqrt{K}}\bhW \bV$ with high probability has singular values close to $1$. 
		
		\begin{proposition} \label{thm:singWV}
			Given $\mathcal{S}$. The singular values of $\frac{1}{\sqrt{K}} \bhW \bV$ with probability at least $1-\delta^*$  lie inside the interval $[\sqrt{1-\varepsilon^*/30},\sqrt{1+\varepsilon^*/30}]$. 
			\begin{proof}
				Define $\tau = \varepsilon^*/30$. 	Notice that $$ K \geq  10^4 l \geq 7.87 \tau^{-2} (6.9 d+\log(r/\delta^*)).$$
				We have, for any $\bx \in \mathbb{R}^d$,
				\begin{equation} \label{eq:singWV1}
				\|\bhW \bV \bx\|_{{2}}^2  
				=\sum^K_{j=1} \|  \sum^r_{i=1} \bhD^{(i,j)} \bH \bD^{(i)} \bV^{(i)} \bx \|_{{2}}^2. 
				\end{equation}
				Denote by $\bd^{(k,j)}$ a vector with $i$-th entry equal to the $(k,k)$-th entry of matrix $\bhD^{(i,j)}$, $1 \leq k \leq r$. Denote by $\bZ^{(k)}$ the matrix with $i$-th row equal to the $k$-th row of matrix $\bH \bD^{(i)} \bV^{(i)}$, $1 \leq k \leq r$. Notice the following relations:
				\begin{equation} \label{eq:singWV2}
				\|  \sum^r_{i=1} \bhD^{(i,j)} \bH \bD^{(i)} \bV^{(i)} \bx \|_{{2}}^2 =  \sum^{r}_{k=1} \langle \bd^{(k,j)},\bZ^{(k)} \bx \rangle^2,~1 \leq j \leq K, 
				\end{equation} 
				and 
				\begin{equation} \label{eq:singWV3} 		
				\sum^{r}_{k=1} \|\bZ^{(k)} \bx \|_{{2}}^2 =  \| \bV \bx \|_{{2}}^2.
				\end{equation}

				We have $\frac{1}{K} \sum^K_{j=1}  \langle  \bd^{(k,j)},\bZ^{(k)}  \bx \rangle^2 = \| \bTheta \bZ^{(k)} \bx \|_{{2}}^2 $, where $\bTheta$ is a $K \times l$ rescaled Rademacher matrix. By~\cite[Proposition 3.7]{balabanov2019randomized}, $\bTheta$ is an $(\tau,\delta^*/r,d)$ OSE, which implies that
				\begin{equation*} 
				\begin{split}
				\forall \bx \in \mathbb{R}^d,~ \frac{1}{K} \sum^K_{j=1}  \langle  \bd^{(k,j)},\bZ^{(k)} \bx \rangle^2  = (1\pm\tau) \| \bZ^{(k)} \bx \|_{{2}}^2,
				\end{split} 
				\end{equation*}
				holds  with probability at least $1-\delta^*/r$. By the summation over $k$ and the union bound argument we conclude that 
				\begin{equation} \label{eq:singWV4}
				\begin{split}
				\forall \bx \in \mathbb{R}^d,~   \frac{1}{K} \sum^K_{j=1} \sum^{r}_{k=1} \langle  \bd^{(k,j)},\bZ^{(k)} \bx \rangle^2= (1 \pm\tau) \sum^{r}_{k=1} \| \bZ^{(k)} \bx \|_{{2}}^2, 
				\end{split}
				\end{equation}
				holds  with probability at least $1-\delta^*$.
				By straightforward substitution of the expressions~\cref{eq:singWV2,eq:singWV3} into~\cref{eq:singWV4}, and using~\cref{eq:singWV1}, we conclude that with probability at least $1-\delta^*$,
				$$\forall \bx \in \mathbb{R}^d,~ \frac{1}{K} \|\bhW \bV \bx\|_{{2}}^2 = (1 \pm\tau) \|\bV \bx\|_{{2}}^2, $$
				which is equivalent to the statement of the proposition.
			\end{proof}    
		\end{proposition}

		\Cref{thm:chernoff} presents a corollary of the Matrix Chernoff inequality from~\cite{tropp2011improved}, used to show that  $\bM = \frac{1}{\sqrt{K}} \bhW \bV$ and $\sqrt{\frac{rK}{{l}}} \bhR \bM = \bOmega \bV$ have approximately equal maximal and minimal singular values.
		
		\begin{proposition}[Corollary of Theorem 2.2 in~\cite{tropp2011improved}] \label{thm:chernoff}
			Let $\bM$ be some $rK \times d$ matrix. Let $0 <\varepsilon^* <1$ and $0<\delta^* <1$. Let {$\bm^{(j)}$} {denote} the rows of $\bM$ and let $M:=rK  \max_{j=1,\hdots,rK} \| {\bm^{(j)}} \|_2^2$ and $N \geq  \sigma_{min}(\bM)^{-2}$. Draw at random a sampling matrix $\bhR$ in~\cref{eq:hRsampling} with
			$$\small {l} \geq 2 ({\varepsilon^*}^{2}-{\varepsilon^*}^{3}/3)^{-1} M N \log (d/\delta^*).$$ Given $\mathcal{S}$, then with probability at least $1-2\delta^*$,
			\small 
			\begin{equation} \label{eq:chernoff}
			\begin{split}
			\sqrt{1-{\varepsilon^*}}\sigma_{min}(\bM) \leq \sigma_{min}(\sqrt{\frac{rK}{{l}}} \bhR \bM )\leq \sigma_{max}(\sqrt{\frac{rK}{{l}}} \bhR \bM ) \leq \sqrt{1+{\varepsilon^*}} \sigma_{max}(\bM). 
			\end{split}
			\end{equation}
			\normalsize
			\begin{proof}
				For any symmetric matrix $\bX$,  let $\lambda_{min}(\bX)$ and $\lambda_{max}(\bX)$ denote the minimal and the maximal eigenvalues of $\bX$.
				To prove~\cref{thm:chernoff} we will use the matrix Chernoff tail bounds from~\cite{tropp2011improved} presented in~\cref{thm:matrixchernoff}.  
				
				Define $X := \{ {\bm^{(j)}} {(\bm^{(j)})^\mathrm{T}} \}^n_{j=1} $. Consider the matrix 
				$$ \bX := (\bhR \bM)^\mathrm{T} \bhR \bM = \sum_{j \in T}{\bm^{(j)}} {(\bm^{(j)})^\mathrm{T}} ,$$
				where $T$ is a set, with $\#T={l}$, of elements of $\{1,2, \hdots, rK \}$ drawn uniformly and without replacement.  The matrix $\bX$ can be written as
				$ \bX = \sum^{{l}}_{i=1} \bX_i,$
				where $\{ \bX_i \}^{{l}}_{i=1}$ is a uniformly drawn, without replacement, random subset of $X$. We have $\mathbb{E}(\bX_1)=\frac{1}{rK} \bM^\mathrm{T} \bM$. Furthermore, \small $$ \lambda_{\mathrm{max}}({\bm^{(j)}} {(\bm^{(j)})^\mathrm{T}}) = \|{\bm^{(j)}} \|^2 \leq \frac{M}{rK},~1 \leq j \leq rK.$$ \normalsize By applying~\Cref{thm:matrixchernoff} and some algebraic operations, we obtain  
				\small 
				\begin{align*}
				\mathbb{P} ({\textstyle\lambda_\mathrm{min}(\bX)}&{\textstyle  \leq (1-{\varepsilon^*}) \lambda_{min}(\bM^\mathrm{T} \bM)  \frac{{l}}{rK}})
				\leq {\textstyle d \left ( \frac{e^{-{\varepsilon^*}}}{(1-{\varepsilon^*})^{1-{\varepsilon^*}}} \right )^{\lambda_{min}(\bM^\mathrm{T} \bM) {l}/M}} \\
				&\leq d~ e^{-({\varepsilon^*}^{2}/2 - {\varepsilon^*}^{3}/6)(MN)^{-1} {l} } \leq \delta,
				\end{align*}
				\begin{align*}
				\mathbb{P} ({\textstyle\lambda_\mathrm{max}(\bX)}&{\textstyle \geq (1+{\varepsilon^*}) \lambda_{max}(\bM^\mathrm{T} \bM)\frac{{l}}{rK}})  \leq {\textstyle d \left ( \frac{e^{{\varepsilon^*}}}{(1+{\varepsilon^*})^{1+{\varepsilon^*}}} \right )^{\lambda_{max}(\bM^\mathrm{T} \bM){l}/M}} \\
				& \leq d~ e^{-({\varepsilon^*}^{2}/2 - {\varepsilon^*}^{3}/6) (MN)^{-1} {l}} \leq \delta.
				\end{align*}
				\normalsize
				The statement of the lemma follows by a union bound 	argument.
			\end{proof}
		\end{proposition}
		
		\begin{theorem}[Matrix Chernoff tail bounds from~\cite{tropp2011improved}] \label{thm:matrixchernoff}
			Consider a finite set $X$ of symmetric positive semi-definite matrices of size $d \times d$.  Define the constant $L:=\max_{\bX_j \in X} \lambda_{max}(\bX_j)$. Let $\{\bX_i\}^{{l}}_{i=1}$ be a uniformly sampled, without replacement, random subset of $X$ and $\bX:=\sum^{{l}}_{i=1}\bX_i$. Then 
			\begin{align*}
			&\mathbb{P}\left ({\textstyle\lambda_{min}(\bX) \leq (1-\varepsilon)\mu_{min} }\right )\leq {\textstyle d \left (\frac{e^{-\varepsilon}}{(1-\varepsilon)^{1-\varepsilon}} \right )^{\mu_{min}/L}} \\
			&\mathbb{P}\left ({\textstyle\lambda_{max}(\bX) \geq (1+\varepsilon)\mu_{max} }\right )\leq {\textstyle d \left (\frac{e^{\varepsilon}}{(1+\varepsilon)^{1+\varepsilon}} \right )^{\mu_{max}/L}}
			\end{align*}
			where $\mu_{min}={l} \lambda_{min}(\mathbb{E}\bX_1)$ and $\mu_{max}= {l} \lambda_{max}(\mathbb{E}\bX_1)$.
		\end{theorem}

		By plugging~\cref{thm:lemma1} and the result of~\cref{thm:singWV} into~\cref{thm:chernoff} and taking $\bM = \frac{1}{\sqrt{K}} \bhW \bV$, $M = (\sqrt{d}+\sqrt{8\log{rK/\delta^*}})^2$, $N =1.07$, along with the union bound argument, we deduce that,  
		\begin{equation*} 
		\begin{split}
		\sqrt{1-{\varepsilon^*}} \sqrt{1-{\varepsilon^*/30}}   
		\leq \sigma_{min}(\bOmega \bV) 
		\leq \sigma_{max}(\bOmega \bV ) 
		\leq \sqrt{1+{\varepsilon^*}} \sqrt{1+{\varepsilon^*/30}}
		\end{split}
		\end{equation*}
		\normalsize
		holds with probability at least $1-4\delta^*$ under the condition $\mathcal{S}$.  	
		Finally by few algebraic operations, we conclude that, given $\mathcal{S}$, the singular values of $\bOmega \bV$ belong to $[\sqrt{1-\varepsilon},\sqrt{1+\varepsilon}]$ with probability at least $4\delta/5$.
		The proof of the main theorem is finished by reminding that $\mathcal{S}$ occurs with probability at least $1-\delta^*$, the union bound argument and few additional algebraic operations.  \qed

		\section{Conclusion} \label{conclusion}
		The proposed block SRHT can combine the advantages of structured and unstructured matrices, such as low application complexity and suitability for distributed computing.
		It should outperform all known embeddings in a distributed architecture with not too large number of processors. At the same time it yields the same approximation guarantees as standard SRHT. 
		We have chosen the low-rank approximation problem as a representative application.  We revised popular randomized methods for this problem with implementation aspects on distributed architectures, and then presented their quasi-optimality characterizations from a {projection-based} point of view, compatible with block SRHT. Numerical validation of the methodology showed that the block SRHT in practice provides solutions of the same quality as Gaussian embeddings. Yet, the block SRHT was up to a factor of $2.5$ faster to apply. Moreover, even greater gains in runtime are expected for larger problems and sampling dimensions.

	    \section{Acknowledgments}
	    This project has received funding from the European Research Council (ERC) under the European Union's Horizon 2020 research and innovation program (grant agreement No 810367). 
		\bibliographystyle{unsrt}
		\small 
		\bibliography{references}
		\normalsize
		
	\end{document}

%% file: accuracy_pot16_srht.tex
          \begin{tikzpicture}
              \datavisualization [
              scientific axes=clean,
              visualize as line=a,
              visualize as line=b,
              visualize as line=c,
              all axes={grid},
              style sheet=strong colors,
              a={style={mark=+}},
              b={style={mark=x}},
              c={style={mark=x}},
              x axis={length=1.8cm,label=Approximation rank,
                     ticks and grid={major={at={200,400,600}}}},
              y axis={logarithmic, include value={1e-5,1e-3},
                      ticks={style={/pgf/number format/retain unit mantissa=false}},
                      label=Trace relative error},
              ]
              data[set=a]{x, y
                200, 0.0004564099468970148
                300, 0.0002674919746696206
                400, 0.0001845297363573884
                500, 0.00015872984766719764
                600, 0.00015324819918411574
              }
              data[set=b]{x, y
                200, 0.0004215880046974093
                300, 0.0002183035081589341
                400, 0.00012125062194734365
                500, 8.469318321043548e-5
                600, 7.297267327279283e-5
              }
              data[set=c]{x, y
                200, 0.00041167079444487197
                300, 0.00020371088805645347
                400, 0.00010218733212497066
                500, 6.163691980252096e-5
                600, 4.692061911134303e-5
              };
            \end{tikzpicture}

%% file: accuracy_pot16_year_1e4.tex
\begin{tikzpicture}
    \datavisualization [
    scientific axes=clean,
    visualize as line=a,
    visualize as line=b,
    visualize as line=c,
    visualize as line=d,
    visualize as line=e,
    all axes={grid},
    style sheet=strong colors,
    a={style={mark=+}},
    b={style={mark=x}},
    c={style={mark=x}},
    d={style={mark=x}},
    e={style={mark=x}},
    legend=east outside,
    x axis={length=3cm,label=Approximation rank,
           ticks and grid={major={at={200,400,600,800,1000}}}},
    y axis={logarithmic, include value={1e-3,1e-2},
            ticks={style={/pgf/number format/retain unit mantissa=false}}},
    ]
    data[set=a]{x, y
      100, 0.00825129036553977
      200, 0.005557963107956253
      300, 0.004568564417985671
      400, 0.004077579841838208
      500, 0.00381757548958801
      600, 0.003681125728413008
    }
    data[set=b]{x, y
      100, 0.007879984845050543
      200, 0.004921808345827181
      300, 0.003738546033823656
      400, 0.0030937007917051828
      500, 0.0027097993330895804
      600, 0.002473215543003254
      700, 0.0023116841929074956
      800, 0.0022035284243754017
      900, 0.0021311785060321676
      1000, 0.0020853506146168944
    }
    data[set=c]{x, y
      100, 0.007718290063396716
      200, 0.004620462899962847
      300, 0.003316855924578431
      400, 0.0025690156548703663
      500, 0.0020947813847318085
      600, 0.0017717339322189769
      700, 0.0015405940191652858
      800, 0.0013704153907022917
      900, 0.0012420652787612162
      1000, 0.0011440733484917765
    }
    data[set=d]{x, y
      100, 0.007701120166296103
      200, 0.00458597367987795
      300, 0.0032664443144641768
      400, 0.002503490839811249
      500, 0.0020154247426186723
      600, 0.0016792321240146834
      700, 0.0014366587751849045
      800, 0.0012555052660278448
      900, 0.0011175014470260634
      1000, 0.0010095411976450389
    }
    data[set=e]{x, y
      100, 0.0076922601988711065
      200, 0.0045690946002515595
      300, 0.003241329367401083
      400, 0.0024698996661912547
      500, 0.0019742536178723076
      600, 0.0016312567186267204
      700, 0.00138199483019422
      800, 0.001194704225178804
      900, 0.0010508637675016638
      1000, 0.0009371617208563574
    };
  \end{tikzpicture}

%% file: accuracy_pot16_year_1e5.tex
\begin{tikzpicture}
      \datavisualization [
      scientific axes=clean,
      visualize as line=a,
      visualize as line=b,
      visualize as line=c,
      visualize as line=d,
      visualize as line=e,
      a={style={mark=+}},
      b={ style={mark=x}},
      c={ style={mark=x}},
      all axes={grid},
      style sheet=strong colors,
      a={label in legend={text={$l = 600$}}, style={mark=+}},
      b={label in legend={text={$l = 1000$}}, style={mark=x}},
      c={label in legend={text={$l = 2000$}}, style={mark=x}},
      d={label in legend={text={$l = 2500$}}, style={mark=x}},
      e={label in legend={text={$l = 3000$}}, style={mark=x}},
      legend=east outside,
      x axis={length=3cm, label=Approximation rank},
      y axis={logarithmic, include value={1e-6,1e-8},
              ticks={style={/pgf/number format/retain unit mantissa=false}}},
      ]
      data[set=a]{x, y
        100, 1.840771566165316e-6
        200, 4.92539556233576e-7
        300, 2.991378097455335e-7
        400, 2.3374390894179435e-7
        500, 2.0626118957345461e-7
        600, 1.9462445629484666e-7
      }
      data[set=b]{x, y
        100, 1.8147141845578402e-6
        200, 4.4405028775796266e-7
        300, 2.3232109523613442e-7
        400, 1.5300754513267565e-7
        500, 1.1464346676972303e-7
        600, 9.333560991536259e-8
        700, 8.124934742716097e-8
        800, 7.403922786392929e-8
        900, 6.974510557479654e-8
        1000, 6.767608335903351e-8
      }
      data[set=c]{x, y
        100, 1.8085236483408485e-6
        200, 4.3176769677952735e-7
        300, 2.1432524069608647e-7
        400, 1.2989735706448955e-7
        500, 8.694000371272181e-8
        700, 4.6360916331071497e-8
        800, 3.60309564240876e-8
        900, 2.890920894197454e-8
        1000, 2.3880767573001803e-8
      }
      data[set=d]{x, y
        100, 1.8081941891245874e-6
        200, 4.3110502140334723e-7
        300, 2.1333659565485474e-7
        400, 1.285883429180018e-7
        500, 8.532278258169787e-8
        600, 6.004036338013848e-8
        700, 4.412755205423226e-8
        800, 3.351877572004217e-8
        900, 2.6133669027101353e-8
        1000, 2.0857336399068046e-8
      }
      data[set=e]{x, y
        100, 1.8080747504025285e-6
        200, 4.3087392652936515e-7
        300, 2.1298960228910387e-7
        400, 1.2812131067986263e-7
        500, 8.473935010178798e-8
        600, 5.934919495204126e-8
        700, 4.332727533264087e-8
        800, 3.260819389215088e-8
        900, 2.5123897532464132e-8
      };

    \end{tikzpicture}

%% file: runtime_julia.tex
\begin{tikzpicture}
    \datavisualization [
    scientific axes=clean,
    visualize as line=a,
    visualize as line=b,
    all axes={grid},
    legend=east outside,
    style sheet=strong colors,
    a={label in legend={text={Gauss.}}, style={mark=x}},
    b={label in legend={text={BSRHT}}, style={mark=x}},
    x axis={length=3.5cm,label=Sampling size,
    ticks and grid={major={at={600,1000,2000}}}},
    y axis={include value=0, label={Runtime (s)}}
    ]
    data[set=b]{x, y
      600 , 2.83
      1000, 3.2
      2000, 4.830
    }
    data[set=a]{x, y
      600 , 3.922
      1000, 5.991
      2000, 12.720
    };
  \end{tikzpicture}

%% file: B.tex
\begin{tikzpicture}
    \datavisualization [
    scientific axes=clean,
    visualize as line=a,
    visualize as line=b,
    visualize as line=c,
    visualize as line=d,
    all axes={grid},
    legend=east outside,
    style sheet=cust,
    a={label in legend={text={Gauss. total}}, style={mark=o}},
    b={label in legend={text={BSRHT total}}, style={mark=o}},
    c={label in legend={text={Gauss. local}}, style={mark=triangle*}},
    d={label in legend={text={BSRHT local}}, style={mark=triangle*}},
    x axis={logarithmic, length=3.5cm,label=\# MPI process,
    ticks and grid={major={at={12,48,192,1536}}}},
    y axis={logarithmic, include value={16,0.055}, label={Runtime (s)}}
    ]
    data[set=d]{x, y
      12, 4.769764E+00
      24, 2.294961E+00
      48, 1.244005E+00
      96, 5.488213E-01
      192, 2.750293E-01
      384, 1.553763E-01
      768, 1.065320E-01
      1536, 5.291957E-02
    }
    data[set=b]{x, y
      12, 5.12066663
      24, 2.4949097
      48, 1.4636852
      96, 0.59322252
      192, 0.317195
      384, 0.19818827
      768, 0.2961686
      1536, 0.22341827
    }
    data[set=c]{x, y
      12, 1.316271E+01
      24, 7.909814E+00
      48, 3.770663E+00
      96, 1.386485E+00
      192, 7.830114E-01
      384, 4.759886E-01
      768, 3.806720E-01
      1536, 2.252002E-01
    }
    data[set=a]{x, y
      12, 15.490014
      24, 10.51908
      48, 4.939651
      96, 1.47704246
      192, 0.82027688
      384, 0.57367648
      768, 0.4993317
      1536, 0.3329018
    };
  \end{tikzpicture}

%% file: C.tex
\begin{tikzpicture}
    \datavisualization [
    scientific axes=clean,
    visualize as line=a,
    visualize as line=b,
    visualize as line=c,
    visualize as line=d,
    all axes={grid},
    legend=east outside,
    style sheet=cust,
    a={label in legend={text={Gauss. total}}, style={mark=o, red}},
    b={label in legend={text={BSRHT total}}, style={mark=o,black}},
    c={label in legend={text={Gauss. local}}, style={mark=triangle*,red}},
    d={label in legend={text={BSRHT local}}, style={mark=triangle*,black}},
    x axis={logarithmic, length=3.5cm,label=\# MPI process,
    ticks and grid={major={at={12,48,192,1536}}}},
    y axis={logarithmic, include value={45,0.1}, label={Runtime (s)}}
    ]
    data[set=d]{x, y
      12, 4.138710E+01
      24, 2.005116E+01
      48, 1.014917E+01
      96, 5.057663E+00
      192, 2.452943E+00
      384, 1.261642E+00
      768, 8.092221E-01
      1536, 4.041584E-01
    }
    data[set=b]{x, y
      12, 44.018897
      24, 21.622836
      48, 11.370489
      96, 5.3587835
      192, 2.54963519
      384, 1.35263526
      768, 1.0106163
      1536, 0.692183
    }
    data[set=c]{x, y
      768, 3.631734E+00
      1536, 1.759755E+00
    }
    data[set=a]{x, y
      768, 4.4825694
      1536, 2.109484
    };
  \end{tikzpicture}

%% file: 100E6_Mem.tex
\begin{tikzpicture}
    \datavisualization [
    scientific axes=clean,
    visualize as line=a,
    visualize as line=b,
    visualize as line=c,
    visualize as line=d,
    all axes={grid},
    legend=east outside,
    style sheet=cust,
    a={label in legend={text={Gauss.}}, style={mark=o}},
    b={label in legend={text={BSRHT}}, style={mark=o}},
    x axis={logarithmic, length=3.5cm,label=\# MPI process,
    ticks and grid={major={at={12,48,192,1536}}}},
    y axis={logarithmic, include value={1E+01, 21E+02}, label={Memory per process (MB)}, ticks={style={/pgf/number format/retain unit mantissa=false}}}
    ]
    data[set=b]{x, y
      12, 13.76E+02
      24, 6.88E+02
      48, 3.44E+02
      96, 1.72E+02
      192, 8.6E+01
      384, 4.4E+01
      768, 2.2E+01
      1536, 1.1E+01
    }
    data[set=a]{x, y
      768, 20.84E+02
      1536, 10.42E+02
    };
  \end{tikzpicture}

%% file: A.tex
\begin{tikzpicture}
    \datavisualization [
    scientific axes=clean,
    visualize as line=a,
    visualize as line=b,
    visualize as line=c,
    visualize as line=d,
    all axes={grid},
    legend=east outside,
    style sheet=cust,
    a={label in legend={text={Gauss. total}}, style={mark=o}},
    b={label in legend={text={BSRHT total}}, style={mark=o}},
    c={label in legend={text={Gauss. local}}, style={mark=triangle*}},
    d={label in legend={text={BSRHT local}}, style={mark=triangle*}},
    x axis={logarithmic, length=3.5cm,label=\# MPI process,
    ticks and grid={major={at={12,48,192,1536}}}},
    y axis={ include value={0.05,1,4},
                      ticks={style={/pgf/number format/retain unit mantissa=false}},
                      label= Runtime (s)},
    ]
    data[set=d]{x, y
      12, 5.503248E-01
      24, 5.671991E-01
      48, 5.442043E-01
      96, 5.378520E-01
      192, 5.684005E-01
      384, 5.960327E-01
      768, 7.386617E-01
      1536, 7.748135E-01
    }
    data[set=b]{x, y
      12, 0.64380597
      24, 0.6770502
      48, 0.6645144
      96, 0.57637693
      192, 0.61613224
      384, 1.1565506
      768, 0.8390045
      1536, 0.9147793
    }
    data[set=c]{x, y
      12, 1.325457E+00
      24, 1.319118E+00
      48, 1.323661E+00
      96, 1.348338E+00
      192, 1.429731E+00
      384, 1.665711E+00
      768, 2.591757E+00
      1536, 2.923965E+00
    }
    data[set=a]{x, y
      12, 1.4622574E+00
      24, 1.4632314E+00
      48, 1.4378729E+00
      96, 1.42978092E+00
      192, 1.4935034E+00
      384, 1.8163747E+00
      768, 3.0361698E+00
      1536, 3.6662905E+00
    };
  \end{tikzpicture}

%% file: Balabanov_Beaupere_Grigori_Lederer_blockSRHT.bbl
\begin{thebibliography}{10}

\bibitem{martinsson2020randomized}
Per-Gunnar Martinsson and Joel~A Tropp.
\newblock Randomized numerical linear algebra: Foundations and algorithms.
\newblock {\em Acta Numerica}, 29:403--572, 2020.

\bibitem{woodruff2014sketching}
David~P Woodruff et~al.
\newblock Sketching as a tool for numerical linear algebra.
\newblock {\em Foundations and Trends{\textregistered} in Theoretical Computer
  Science}, 10(1--2):1--157, 2014.

\bibitem{vershynin2018high}
Roman Vershynin.
\newblock {\em High-dimensional probability: An introduction with applications
  in data science}, volume~47.
\newblock Cambridge university press, 2018.

\bibitem{mahoney2011randomized}
Michael~W Mahoney.
\newblock Randomized algorithms for matrices and data.
\newblock {\em arXiv preprint arXiv:1104.5557}, 2011.

\bibitem{johnson1984extensions}
William~B Johnson and Joram Lindenstrauss.
\newblock Extensions of lipschitz mappings into a hilbert space 26.
\newblock {\em Contemporary mathematics}, 26:28, 1984.

\bibitem{ailon2006approximate}
Nir Ailon and Bernard Chazelle.
\newblock Approximate nearest neighbors and the fast johnson-lindenstrauss
  transform.
\newblock In {\em Proceedings of the thirty-eighth annual ACM symposium on
  Theory of computing}, pages 557--563, 2006.

\bibitem{sarlos2006improved}
Tamas Sarlos.
\newblock Improved approximation algorithms for large matrices via random
  projections.
\newblock In {\em 2006 47th annual IEEE symposium on foundations of computer
  science (FOCS'06)}, pages 143--152. IEEE, 2006.

\bibitem{halko2011finding}
Nathan Halko, Per-Gunnar Martinsson, and Joel~A Tropp.
\newblock Finding structure with randomness: Probabilistic algorithms for
  constructing approximate matrix decompositions.
\newblock {\em SIAM review}, 53(2):217--288, 2011.

\bibitem{sun2020low}
Yiming Sun, Yang Guo, Charlene Luo, Joel Tropp, and Madeleine Udell.
\newblock Low-rank tucker approximation of a tensor from streaming data.
\newblock {\em SIAM Journal on Mathematics of Data Science}, 2(4):1123--1150,
  2020.

\bibitem{anaraki2013compressive}
Farhad~Pourkamali Anaraki and Shannon~M Hughes.
\newblock Compressive k-svd.
\newblock In {\em 2013 IEEE International Conference on Acoustics, Speech and
  Signal Processing}, pages 5469--5473. IEEE, 2013.

\bibitem{balabanov2019randomized}
Oleg Balabanov and Anthony Nouy.
\newblock Randomized linear algebra for model reduction. {Part} {I}: {Galerkin}
  methods and error estimation.
\newblock {\em Advances in Computational Mathematics}, 45(5-6):2969--3019,
  December 2019.

\bibitem{balabanov2020randomized}
Oleg Balabanov and Laura Grigori.
\newblock Randomized gram-schmidt process with application to gmres.
\newblock {\em arXiv preprint arXiv:2011.05090}, 2020.

\bibitem{alaoui2015fast}
Ahmed Alaoui and Michael~W Mahoney.
\newblock Fast randomized kernel ridge regression with statistical guarantees.
\newblock {\em Advances in neural information processing systems}, 28, 2015.

\bibitem{bach2013sharp}
Francis Bach.
\newblock Sharp analysis of low-rank kernel matrix approximations.
\newblock In {\em Conference on Learning Theory}, pages 185--209. PMLR, 2013.

\bibitem{derezinski2020improved}
Michal Derezinski, Rajiv Khanna, and Michael~W Mahoney.
\newblock Improved guarantees and a multiple-descent curve for column subset
  selection and the nystrom method.
\newblock {\em Advances in Neural Information Processing Systems},
  33:4953--4964, 2020.

\bibitem{meanti2020kernel}
Giacomo Meanti, Luigi Carratino, Lorenzo Rosasco, and Alessandro Rudi.
\newblock Kernel methods through the roof: handling billions of points
  efficiently.
\newblock {\em Advances in Neural Information Processing Systems},
  33:14410--14422, 2020.

\bibitem{yin2021distributed}
Rong Yin, Weiping Wang, and Dan Meng.
\newblock Distributed nystr{\"o}m kernel learning with communications.
\newblock In {\em International Conference on Machine Learning}, pages
  12019--12028. PMLR, 2021.

\bibitem{rudi2017falkon}
Alessandro Rudi, Luigi Carratino, and Lorenzo Rosasco.
\newblock Falkon: An optimal large scale kernel method.
\newblock {\em Advances in neural information processing systems}, 30, 2017.

\bibitem{zhang2013divide}
Yuchen Zhang, John Duchi, and Martin Wainwright.
\newblock Divide and conquer kernel ridge regression.
\newblock In {\em Conference on learning theory}, pages 592--617. PMLR, 2013.

\bibitem{calandriello2016analysis}
Daniele Calandriello, Alessandro Lazaric, and Michal Valko.
\newblock Analysis of nystr{\"o}m method with sequential ridge leverage score
  sampling.
\newblock 2016.

\bibitem{tropp2017fixed}
Joel~A Tropp, Alp Yurtsever, Madeleine Udell, and Volkan Cevher.
\newblock Fixed-rank approximation of a positive-semidefinite matrix from
  streaming data.
\newblock {\em Advances in Neural Information Processing Systems}, 30, 2017.

\bibitem{gittens2013revisiting}
Alex Gittens and Michael Mahoney.
\newblock Revisiting the nystrom method for improved large-scale machine
  learning.
\newblock In {\em International Conference on Machine Learning}, pages
  567--575. PMLR, 2013.

\bibitem{upadhyay2016fast}
Jalaj Upadhyay.
\newblock Fast and space-optimal low-rank factorization in the streaming model
  with application in differential privacy.
\newblock {\em arXiv preprint arXiv:1604.01429}, 2016.

\bibitem{tropp2019streaming}
Joel~A Tropp, Alp Yurtsever, Madeleine Udell, and Volkan Cevher.
\newblock Streaming low-rank matrix approximation with an application to
  scientific simulation.
\newblock {\em SIAM Journal on Scientific Computing}, 41(4):A2430--A2463, 2019.

\bibitem{kannan2014principal}
Ravi Kannan, Santosh Vempala, and David Woodruff.
\newblock Principal component analysis and higher correlations for distributed
  data.
\newblock In {\em Conference on Learning Theory}, pages 1040--1057. PMLR, 2014.

\bibitem{tropp2011improved}
Joel~A. Tropp.
\newblock Improved analysis of the subsampled randomized {Hadamard} transform.
\newblock {\em arXiv:1011.1595 [cs, math]}, July 2011.
\newblock arXiv: 1011.1595.

\bibitem{boutsidis2013improved}
Christos Boutsidis and Alex Gittens.
\newblock Improved matrix algorithms via the subsampled randomized hadamard
  transform.
\newblock {\em SIAM Journal on Matrix Analysis and Applications},
  34(3):1301--1340, 2013.

\bibitem{yang2015implementing}
Jiyan Yang, Xiangrui Meng, and Michael~W Mahoney.
\newblock Implementing randomized matrix algorithms in parallel and distributed
  environments.
\newblock {\em Proceedings of the IEEE}, 104(1):58--92, 2015.

\bibitem{demmel2012communication}
James Demmel, Laura Grigori, Mark Hoemmen, and Julien Langou.
\newblock Communication-optimal parallel and sequential qr and lu
  factorizations.
\newblock {\em SIAM Journal on Scientific Computing}, 34(1):A206--A239, 2012.

\bibitem{chiu2013sublinear}
Jiawei Chiu and Laurent Demanet.
\newblock Sublinear randomized algorithms for skeleton decompositions.
\newblock {\em SIAM Journal on Matrix Analysis and Applications},
  34(3):1361--1383, 2013.

\bibitem{drineas2005nystrom}
Petros Drineas, Michael~W Mahoney, and Nello Cristianini.
\newblock On the nystr{\"o}m method for approximating a gram matrix for
  improved kernel-based learning.
\newblock {\em journal of machine learning research}, 6(12), 2005.

\bibitem{pourkamali2019improved}
Farhad Pourkamali-Anaraki and Stephen Becker.
\newblock Improved fixed-rank nystr{\"o}m approximation via qr decomposition:
  Practical and theoretical aspects.
\newblock {\em Neurocomputing}, 363:261--272, 2019.

\bibitem{pourkamali2018randomized}
Farhad Pourkamali-Anaraki, Stephen Becker, and Michael Wakin.
\newblock Randomized clustered nystrom for large-scale kernel machines.
\newblock In {\em Proceedings of the AAAI Conference on Artificial
  Intelligence}, volume~32, 2018.

\bibitem{li2017algorithm}
Huamin Li, George~C Linderman, Arthur Szlam, Kelly~P Stanton, Yuval Kluger, and
  Mark Tygert.
\newblock Algorithm 971: An implementation of a randomized algorithm for
  principal component analysis.
\newblock {\em ACM Transactions on Mathematical Software (TOMS)}, 43(3):1--14,
  2017.

\bibitem{gittens2011spectral}
Alex Gittens.
\newblock The spectral norm error of the naive nystrom extension.
\newblock {\em arXiv preprint arXiv:1110.5305}, 2011.

\bibitem{tropp2017practical}
Joel~A Tropp, Alp Yurtsever, Madeleine Udell, and Volkan Cevher.
\newblock Practical sketching algorithms for low-rank matrix approximation.
\newblock {\em SIAM Journal on Matrix Analysis and Applications},
  38(4):1454--1485, 2017.

\bibitem{726791}
Y.~Lecun, L.~Bottou, Y.~Bengio, and P.~Haffner.
\newblock Gradient-based learning applied to document recognition.
\newblock {\em Proceedings of the IEEE}, 86(11):2278--2324, 1998.

\bibitem{Bertin-Mahieux2011}
Thierry Bertin-Mahieux, Daniel~P.W. Ellis, Brian Whitman, and Paul Lamere.
\newblock The million song dataset.
\newblock In {\em {Proceedings of the 12th International Conference on Music
  Information Retrieval ({ISMIR} 2011)}}, 2011.

\end{thebibliography}
